
\mag 1200

\input amstex

\expandafter\ifx\csname beta.def\endcsname\relax \else\endinput\fi
\expandafter\edef\csname beta.def\endcsname{%
 \catcode`\noexpand\@=\the\catcode`\@\space}

\let\atbefore @

\catcode`\@=11

\overfullrule\z@

\def\PaperA4{\hsize 6.25truein \vsize 9.63truein}


\def\foliorm{\ifMag\eightrm\else\ninerm\fi}

\let\@ft@\expandafter \let\@tb@f@\atbefore

\newif\ifMag
\def\Magset{\ifnum\mag>\@m\Magtrue\fi}
\Magset

\newif\ifUS

\newdimen\p@@ \p@@\p@
\def\m@ths@r{\ifnum\mathsurround=\z@\z@\else\maths@r\fi}
\def\maths@r{1.6\p@@} \def\mathsurzero{\def\maths@r{\z@}}

\mathsurround\maths@r
\font\Brm=cmr12 \font\Bbf=cmbx12 \font\Bit=cmti12 \font\ssf=cmss10
\font\Bsl=cmsl10 scaled 1200 \font\Bmmi=cmmi10 scaled 1200
\font\BBf=cmbx12 scaled 1200 \font\BMmi=cmmi10 scaled 1440

\def\atletter{\edef\atrestore{\catcode`\noexpand\@=\the\catcode`\@\space}
 \catcode`\@=11}

\newread\@ux \newwrite\@@x \newwrite\@@cd
\let\@np@@\input
\def\@np@t#1{\openin\@ux#1\relax\ifeof\@ux\else\closein\@ux\relax\@np@@ #1\fi}
\def\input#1 {\openin\@ux#1\relax\ifeof\@ux\wrs@x{No file #1}\else
 \closein\@ux\relax\@np@@ #1\fi}
\def\Input#1 {\relax} 

\def\wr@@x#1{} \def\wrs@x{\immediate\write\sixt@@n}

\def\readldf{\@np@t{\jobname.ldf}}
\def\writeldf{\def\wr@@x{\immediate\write\@@x}\def\wr@x@{\write\@@x}
 \def\cl@selbl{\wr@@x{\string\def\string\nextpage{\the\pageno}}%
 \wr@@x{\string\endinput}\immediate\closeout\@@x}
 \immediate\openout\@@x\jobname.ldf}
\let\cl@selbl\relax

\def\nextpage{1}

\def\tod@y{\ifcase\month\or
 January\or February\or March\or April\or May\or June\or July\or
 August\or September\or October\or November\or December\fi\space\,
\number\day,\space\,\number\year}

\newcount\c@time
\def\h@@r{hh}\def\m@n@te{mm}
\def\wh@tt@me{\c@time\time\divide\c@time 60\edef\h@@r{\number\c@time}%
 \multiply\c@time -60\advance\c@time\time\edef
 \m@n@te{\ifnum\c@time<10 0\fi\number\c@time}}
\def\t@me{\h@@r\/{\rm:}\m@n@te} \let\whattime\wh@tt@me
\def\today{\tod@y\wr@@x{\string\todaydef{\tod@y}}}
\def\nowtime{\t@me{\let\/\ic@\wr@@x{\string\nowtimedef{\t@me}}}}
\def\todaydef#1{} \def\nowtimedef#1{}

\def\em#1{{\it #1\/}} \def\emph#1{{\sl #1\/}}

\def\fitem#1{\par\setbox\z@\hbox{#1}\hangindent\wd\z@
 \hglue-2\parindent\kern\wd\z@\indent\llap{#1}\ignore}

\def\itemflat#1{\par\setbox\z@\hbox{\rm #1\enspace}\hang\ifnum\wd\z@>\parindent
 \noindent\unhbox\z@\ignore\else\textindent{\rm#1}\fi}

\newcount\itemlet
\def\newbi{\itemlet 96} \newbi
\def\bitem{\gad\itemlet \par\hangindent1.5\parindent
 \hglue-.5\parindent\textindent{\rm\rlap{\char\the\itemlet}\hp{b})}}
\def\atem{\newbi\bitem}

\newcount\itemrm

\def\iitem{\gad\itemrm \par\hangindent1.5\parindent
 \hglue-.5\parindent\textindent{\rm\hp{v}\llap{\romannumeral\the\itemrm})}}

\newcount\itemar

\def\iitema{\gad\itemrm \par\hangindent1.5\parindent
 \hglue-.5\parindent\textindent{\rm\hp{0}\llap{\the\itemrm}.}}

\def\center{\par\begingroup\leftskip\z@ plus \hsize \rightskip\leftskip
 \parindent\z@\parfillskip\z@skip \def\\{\unskip\break}}
\def\endcenter{\endgraf\endgroup}

\let\b@gr@@\begingroup \let\B@gr@@\begingroup
\def\b@gr@{\b@gr@@\let\b@gr@@\undefined}
\def\B@gr@{\B@gr@@\let\B@gr@@\undefined}

\def\@fn@xt#1#2#3{\let\@ch@r=#1\def\n@xt{\ifx\t@st@\@ch@r
 \def\n@@xt{#2}\else\def\n@@xt{#3}\fi\n@@xt}\futurelet\t@st@\n@xt}

\def\@fwd@@#1#2#3{\setbox\z@\hbox{#1}\ifdim\wd\z@>\z@#2\else#3\fi}
\def\s@twd@#1#2{\setbox\z@\hbox{#2}#1\wd\z@}

\def\r@st@re#1{\let#1\s@v@} \def\s@v@d@f{\let\s@v@}

\def\p@sk@p#1#2{\par\skip@#2\relax\ifdim\lastskip<\skip@\relax\removelastskip
 \ifnum#1=\z@\else\penalty#1\relax\fi\vskip\skip@
 \else\ifnum#1=\z@\else\penalty#1\relax\fi\fi}
\def\sk@@p#1{\par\skip@#1\relax\ifdim\lastskip<\skip@\relax\removelastskip
 \vskip\skip@\fi}

\newbox\p@b@ld
\def\poorbold#1{\setbox\p@b@ld\hbox{#1}\kern-.01em\copy\p@b@ld\kern-\wd\p@b@ld
 \kern.02em\copy\p@b@ld\kern-\wd\p@b@ld\kern-.012em\raise.02em\box\p@b@ld}

\ifx\plainfootnote\undefined \let\plainfootnote\footnote \fi

\let\s@v@\proclaim \let\proclaim\relax
\def\r@R@fs#1{\let#1\s@R@fs} \let\s@R@fs\Refs \let\Refs\relax
\def\r@endd@#1{\let#1\s@endd@} \let\s@endd@\enddocument
\let\bye\relax

\def\myR@fs{\@fn@xt[\m@R@f@\m@R@fs} \def\m@R@fs{\@fn@xt*\m@r@f@@\m@R@f@@}
\def\m@R@f@@{\m@R@f@[References]} \def\m@r@f@@*{\m@R@f@[]}

\def\Twelvepoint{\twelvepoint \let\Bbf\BBf \let\Bmmi\BMmi
\font\Brm=cmr12 scaled 1200 \font\Bit=cmti12 scaled 1200
\font\ssf=cmss10 scaled 1200 \font\Bsl=cmsl10 scaled 1440
\font\BBf=cmbx12 scaled 1440 \font\BMmi=cmmi10 scaled 1728}

\newdimen\b@gsize

\newdimen\r@f@nd \newbox\r@f@b@x \newbox\adjb@x
\newbox\p@nct@ \newbox\k@yb@x \newcount\rcount
\newbox\b@b@x \newbox\p@p@rb@x \newbox\j@@rb@x \newbox\y@@rb@x
\newbox\v@lb@x \newbox\is@b@x \newbox\p@g@b@x \newif\ifp@g@ \newif\ifp@g@s
\newbox\inb@@kb@x \newbox\b@@kb@x \newbox\p@blb@x \newbox\p@bl@db@x
\newbox\ed@b@x \newif\ifed@ \newif\ifed@s \newif\if@fl@b \newif\if@fn@m
\newbox\p@p@nf@b@x \newbox\inf@b@x \newbox\b@@nf@b@x
\newtoks\@dd@p@n \newtoks\@ddt@ks

\newif\ifp@gen@

\def\p@@nt{.\kern.3em} \let\point\p@@nt

\let\proheadfont\bf \let\probodyfont\sl \let\demofont\it

\headline={\hfil}
\footline={\ifp@gen@\ifnum\pageno=\z@\else\hfil\foliorm\folio\fi\else
 \ifnum\pageno=\z@\hfil\foliorm\folio\fi\fi\hfil\global\p@gen@true}
\parindent1pc

\font@\tensmc=cmcsc10
\font@\sevenex=cmex7
\font@\sevenit=cmti7
\font@\eightrm=cmr8
\font@\sixrm=cmr6
\font@\eighti=cmmi8 \skewchar\eighti='177
\font@\sixi=cmmi6 \skewchar\sixi='177
\font@\eightsy=cmsy8 \skewchar\eightsy='60
\font@\sixsy=cmsy6 \skewchar\sixsy='60
\font@\eightex=cmex8
\font@\eightbf=cmbx8
\font@\sixbf=cmbx6
\font@\eightit=cmti8
\font@\eightsl=cmsl8
\font@\eightsmc=cmcsc8
\font@\eighttt=cmtt8
\font@\ninerm=cmr9
\font@\ninei=cmmi9 \skewchar\ninei='177
\font@\ninesy=cmsy9 \skewchar\ninesy='60
\font@\nineex=cmex9
\font@\ninebf=cmbx9
\font@\nineit=cmti9
\font@\ninesl=cmsl9
\font@\ninesmc=cmcsc9
\font@\ninemsa=msam9
\font@\ninemsb=msbm9
\font@\nineeufm=eufm9
\font@\eightmsa=msam8
\font@\eightmsb=msbm8
\font@\eighteufm=eufm8
\font@\sixmsa=msam6
\font@\sixmsb=msbm6
\font@\sixeufm=eufm6

\loadmsam\loadmsbm\loadeufm
\input amssym.tex

\def\footnoterule{\kern-3\p@\hrule width5pc\kern 2.6\p@}
\def\m@k@foot#1{\insert\footins
 {\interlinepenalty\interfootnotelinepenalty
 \ifMag\eightpoint\else\ninepoint\fi
 \splittopskip\ht\strutbox\splitmaxdepth\dp\strutbox
 \floatingpenalty\@MM\leftskip\z@\rightskip\z@
 \spaceskip\z@\xspaceskip\z@
 \leavevmode\footstrut\ignore#1\unskip\lower\dp\strutbox
 \vbox to\dp\strutbox{}}}
\def\ftext#1{\m@k@foot{\vsk-.8>\nt #1}}
\def\pr@cl@@m#1{\p@sk@p{-100}\medskipamount
 \def\endproclaim{\endgroup\p@sk@p{55}\medskipamount}\begingroup
 \nt\ignore\proheadfont#1\unskip.\enspace\probodyfont\ignore}
\outer\def\proclaim{\pr@cl@@m} \s@v@d@f\proclaim \let\proclaim\relax
\def\demo#1{\sk@@p\medskipamount\nt{\ignore\demofont#1\unskip.}\enspace
 \ignore}
\def\enddemo{\sk@@p\medskipamount}

\def\cite#1{{\rm[#1]}} 
 \def\Refs#1#2{\relax}

\def\big@#1#2{{\hbox{$\left#2\vcenter to#1\b@gsize{}%
 \right.\nulldelimiterspace\z@\m@th$}}}
\def\big{\big@\@ne}
\def\Big{\big@{1.5}}
\def\bigg{\big@\tw@}
\def\Bigg{\big@{2.5}}
\normallineskiplimit\p@

\def\tenpoint{\p@@\p@ \normallineskiplimit\p@@
 \mathsurround\m@ths@r \normalbaselineskip12\p@@
 \abovedisplayskip12\p@@ plus3\p@@ minus9\p@@
 \belowdisplayskip\abovedisplayskip
 \abovedisplayshortskip\z@ plus3\p@@
 \belowdisplayshortskip7\p@@ plus3\p@@ minus4\p@@
 \textonlyfont@\rm\tenrm \textonlyfont@\it\tenit
 \textonlyfont@\sl\tensl \textonlyfont@\bf\tenbf
 \textonlyfont@\smc\tensmc \textonlyfont@\tt\tentt
 \ifsyntax@ \def\big##1{{\hbox{$\left##1\right.$}}}%
  \let\Big\big \let\bigg\big \let\Bigg\big
 \else
  \textfont\z@\tenrm \scriptfont\z@\sevenrm \scriptscriptfont\z@\fiverm
  \textfont\@ne\teni \scriptfont\@ne\seveni \scriptscriptfont\@ne\fivei
  \textfont\tw@\tensy \scriptfont\tw@\sevensy \scriptscriptfont\tw@\fivesy
  \textfont\thr@@\tenex \scriptfont\thr@@\sevenex
	\scriptscriptfont\thr@@\sevenex
  \textfont\itfam\tenit \scriptfont\itfam\sevenit
	\scriptscriptfont\itfam\sevenit
  \textfont\bffam\tenbf \scriptfont\bffam\sevenbf
	\scriptscriptfont\bffam\fivebf
  \textfont\msafam\tenmsa \scriptfont\msafam\sevenmsa
	\scriptscriptfont\msafam\fivemsa
  \textfont\msbfam\tenmsb \scriptfont\msbfam\sevenmsb
	\scriptscriptfont\msbfam\fivemsb
  \textfont\eufmfam\teneufm \scriptfont\eufmfam\seveneufm
	\scriptscriptfont\eufmfam\fiveeufm
  \setbox\strutbox\hbox{\vrule height8.5\p@@ depth3.5\p@@ width\z@}%
  \setbox\strutbox@\hbox{\lower.5\normallineskiplimit\vbox{%
	\kern-\normallineskiplimit\copy\strutbox}}%
   \setbox\z@\vbox{\hbox{$($}\kern\z@}\b@gsize1.2\ht\z@
  \fi
  \normalbaselines\rm\dotsspace@1.5mu\ex@.2326ex\jot3\ex@}

\def\eightpoint{\p@@.8\p@ \normallineskiplimit\p@@
 \mathsurround\m@ths@r \normalbaselineskip10\p@
 \abovedisplayskip10\p@ plus2.4\p@ minus7.2\p@
 \belowdisplayskip\abovedisplayskip
 \abovedisplayshortskip\z@ plus3\p@@
 \belowdisplayshortskip7\p@@ plus3\p@@ minus4\p@@
 \textonlyfont@\rm\eightrm \textonlyfont@\it\eightit
 \textonlyfont@\sl\eightsl \textonlyfont@\bf\eightbf
 \textonlyfont@\smc\eightsmc \textonlyfont@\tt\eighttt
 \ifsyntax@\def\big##1{{\hbox{$\left##1\right.$}}}%
  \let\Big\big \let\bigg\big \let\Bigg\big
 \else
  \textfont\z@\eightrm \scriptfont\z@\sixrm \scriptscriptfont\z@\fiverm
  \textfont\@ne\eighti \scriptfont\@ne\sixi \scriptscriptfont\@ne\fivei
  \textfont\tw@\eightsy \scriptfont\tw@\sixsy \scriptscriptfont\tw@\fivesy
  \textfont\thr@@\eightex \scriptfont\thr@@\sevenex
	\scriptscriptfont\thr@@\sevenex
  \textfont\itfam\eightit \scriptfont\itfam\sevenit
	\scriptscriptfont\itfam\sevenit
  \textfont\bffam\eightbf \scriptfont\bffam\sixbf
	\scriptscriptfont\bffam\fivebf
  \textfont\msafam\eightmsa \scriptfont\msafam\sixmsa
	\scriptscriptfont\msafam\fivemsa
  \textfont\msbfam\eightmsb \scriptfont\msbfam\sixmsb
	\scriptscriptfont\msbfam\fivemsb
  \textfont\eufmfam\eighteufm \scriptfont\eufmfam\sixeufm
	\scriptscriptfont\eufmfam\fiveeufm
 \setbox\strutbox\hbox{\vrule height7\p@ depth3\p@ width\z@}%
 \setbox\strutbox@\hbox{\raise.5\normallineskiplimit\vbox{%
   \kern-\normallineskiplimit\copy\strutbox}}%
 \setbox\z@\vbox{\hbox{$($}\kern\z@}\b@gsize1.2\ht\z@
 \fi
 \normalbaselines\eightrm\dotsspace@1.5mu\ex@.2326ex\jot3\ex@}

\def\ninepoint{\p@@.9\p@ \normallineskiplimit\p@@
 \mathsurround\m@ths@r \normalbaselineskip11\p@
 \abovedisplayskip11\p@ plus2.7\p@ minus8.1\p@
 \belowdisplayskip\abovedisplayskip
 \abovedisplayshortskip\z@ plus3\p@@
 \belowdisplayshortskip7\p@@ plus3\p@@ minus4\p@@
 \textonlyfont@\rm\ninerm \textonlyfont@\it\nineit
 \textonlyfont@\sl\ninesl \textonlyfont@\bf\ninebf
 \textonlyfont@\smc\ninesmc \textonlyfont@\tt\ninett
 \ifsyntax@ \def\big##1{{\hbox{$\left##1\right.$}}}%
  \let\Big\big \let\bigg\big \let\Bigg\big
 \else
  \textfont\z@\ninerm \scriptfont\z@\sevenrm \scriptscriptfont\z@\fiverm
  \textfont\@ne\ninei \scriptfont\@ne\seveni \scriptscriptfont\@ne\fivei
  \textfont\tw@\ninesy \scriptfont\tw@\sevensy \scriptscriptfont\tw@\fivesy
  \textfont\thr@@\nineex \scriptfont\thr@@\sevenex
	\scriptscriptfont\thr@@\sevenex
  \textfont\itfam\nineit \scriptfont\itfam\sevenit
	\scriptscriptfont\itfam\sevenit
  \textfont\bffam\ninebf \scriptfont\bffam\sevenbf
	\scriptscriptfont\bffam\fivebf
  \textfont\msafam\ninemsa \scriptfont\msafam\sevenmsa
	\scriptscriptfont\msafam\fivemsa
  \textfont\msbfam\ninemsb \scriptfont\msbfam\sevenmsb
	\scriptscriptfont\msbfam\fivemsb
  \textfont\eufmfam\nineeufm \scriptfont\eufmfam\seveneufm
	\scriptscriptfont\eufmfam\fiveeufm
  \setbox\strutbox\hbox{\vrule height8.5\p@@ depth3.5\p@@ width\z@}%
  \setbox\strutbox@\hbox{\lower.5\normallineskiplimit\vbox{%
	\kern-\normallineskiplimit\copy\strutbox}}%
   \setbox\z@\vbox{\hbox{$($}\kern\z@}\b@gsize1.2\ht\z@
  \fi
  \normalbaselines\rm\dotsspace@1.5mu\ex@.2326ex\jot3\ex@}

\font@\twelverm=cmr10 scaled 1200
\font@\twelveit=cmti10 scaled 1200
\font@\twelvesl=cmsl10 scaled 1200
\font@\twelvebf=cmbx10 scaled 1200
\font@\twelvesmc=cmcsc10 scaled 1200
\font@\twelvett=cmtt10 scaled 1200
\font@\twelvei=cmmi10 scaled 1200 \skewchar\twelvei='177
\font@\twelvesy=cmsy10 scaled 1200 \skewchar\twelvesy='60
\font@\twelveex=cmex10 scaled 1200
\font@\twelvemsa=msam10 scaled 1200
\font@\twelvemsb=msbm10 scaled 1200
\font@\twelveeufm=eufm10 scaled 1200

\def\twelvepoint{\p@@1.2\p@ \normallineskiplimit\p@@
 \mathsurround\m@ths@r \normalbaselineskip12\p@@
 \abovedisplayskip12\p@@ plus3\p@@ minus9\p@@
 \belowdisplayskip\abovedisplayskip
 \abovedisplayshortskip\z@ plus3\p@@
 \belowdisplayshortskip7\p@@ plus3\p@@ minus4\p@@
 \textonlyfont@\rm\twelverm \textonlyfont@\it\twelveit
 \textonlyfont@\sl\twelvesl \textonlyfont@\bf\twelvebf
 \textonlyfont@\smc\twelvesmc \textonlyfont@\tt\twelvett
 \ifsyntax@ \def\big##1{{\hbox{$\left##1\right.$}}}%
  \let\Big\big \let\bigg\big \let\Bigg\big
 \else
  \textfont\z@\twelverm \scriptfont\z@\eightrm \scriptscriptfont\z@\sixrm
  \textfont\@ne\twelvei \scriptfont\@ne\eighti \scriptscriptfont\@ne\sixi
  \textfont\tw@\twelvesy \scriptfont\tw@\eightsy \scriptscriptfont\tw@\sixsy
  \textfont\thr@@\twelveex \scriptfont\thr@@\eightex
	\scriptscriptfont\thr@@\sevenex
  \textfont\itfam\twelveit \scriptfont\itfam\eightit
	\scriptscriptfont\itfam\sevenit
  \textfont\bffam\twelvebf \scriptfont\bffam\eightbf
	\scriptscriptfont\bffam\sixbf
  \textfont\msafam\twelvemsa \scriptfont\msafam\eightmsa
	\scriptscriptfont\msafam\sixmsa
  \textfont\msbfam\twelvemsb \scriptfont\msbfam\eightmsb
	\scriptscriptfont\msbfam\sixmsb
  \textfont\eufmfam\twelveeufm \scriptfont\eufmfam\eighteufm
	\scriptscriptfont\eufmfam\sixeufm
  \setbox\strutbox\hbox{\vrule height8.5\p@@ depth3.5\p@@ width\z@}%
  \setbox\strutbox@\hbox{\lower.5\normallineskiplimit\vbox{%
	\kern-\normallineskiplimit\copy\strutbox}}%
  \setbox\z@\vbox{\hbox{$($}\kern\z@}\b@gsize1.2\ht\z@
  \fi
  \normalbaselines\rm\dotsspace@1.5mu\ex@.2326ex\jot3\ex@}

\font@\twelvetrm=cmr10 at 12truept
\font@\twelvetit=cmti10 at 12truept
\font@\twelvetsl=cmsl10 at 12truept
\font@\twelvetbf=cmbx10 at 12truept
\font@\twelvetsmc=cmcsc10 at 12truept
\font@\twelvettt=cmtt10 at 12truept
\font@\twelveti=cmmi10 at 12truept \skewchar\twelveti='177
\font@\twelvetsy=cmsy10 at 12truept \skewchar\twelvetsy='60
\font@\twelvetex=cmex10 at 12truept
\font@\twelvetmsa=msam10 at 12truept
\font@\twelvetmsb=msbm10 at 12truept
\font@\twelveteufm=eufm10 at 12truept

\def\twelvetruepoint{\p@@1.2truept \normallineskiplimit\p@@
 \mathsurround\m@ths@r \normalbaselineskip12\p@@
 \abovedisplayskip12\p@@ plus3\p@@ minus9\p@@
 \belowdisplayskip\abovedisplayskip
 \abovedisplayshortskip\z@ plus3\p@@
 \belowdisplayshortskip7\p@@ plus3\p@@ minus4\p@@
 \textonlyfont@\rm\twelvetrm \textonlyfont@\it\twelvetit
 \textonlyfont@\sl\twelvetsl \textonlyfont@\bf\twelvetbf
 \textonlyfont@\smc\twelvetsmc \textonlyfont@\tt\twelvettt
 \ifsyntax@ \def\big##1{{\hbox{$\left##1\right.$}}}%
  \let\Big\big \let\bigg\big \let\Bigg\big
 \else
  \textfont\z@\twelvetrm \scriptfont\z@\eightrm \scriptscriptfont\z@\sixrm
  \textfont\@ne\twelveti \scriptfont\@ne\eighti \scriptscriptfont\@ne\sixi
  \textfont\tw@\twelvetsy \scriptfont\tw@\eightsy \scriptscriptfont\tw@\sixsy
  \textfont\thr@@\twelvetex \scriptfont\thr@@\eightex
	\scriptscriptfont\thr@@\sevenex
  \textfont\itfam\twelvetit \scriptfont\itfam\eightit
	\scriptscriptfont\itfam\sevenit
  \textfont\bffam\twelvetbf \scriptfont\bffam\eightbf
	\scriptscriptfont\bffam\sixbf
  \textfont\msafam\twelvetmsa \scriptfont\msafam\eightmsa
	\scriptscriptfont\msafam\sixmsa
  \textfont\msbfam\twelvetmsb \scriptfont\msbfam\eightmsb
	\scriptscriptfont\msbfam\sixmsb
  \textfont\eufmfam\twelveteufm \scriptfont\eufmfam\eighteufm
	\scriptscriptfont\eufmfam\sixeufm
  \setbox\strutbox\hbox{\vrule height8.5\p@@ depth3.5\p@@ width\z@}%
  \setbox\strutbox@\hbox{\lower.5\normallineskiplimit\vbox{%
	\kern-\normallineskiplimit\copy\strutbox}}%
  \setbox\z@\vbox{\hbox{$($}\kern\z@}\b@gsize1.2\ht\z@
  \fi
  \normalbaselines\rm\dotsspace@1.5mu\ex@.2326ex\jot3\ex@}

\font@\elevenrm=cmr10 scaled 1095
\font@\elevenit=cmti10 scaled 1095
\font@\elevensl=cmsl10 scaled 1095
\font@\elevenbf=cmbx10 scaled 1095
\font@\elevensmc=cmcsc10 scaled 1095
\font@\eleventt=cmtt10 scaled 1095
\font@\eleveni=cmmi10 scaled 1095 \skewchar\eleveni='177
\font@\elevensy=cmsy10 scaled 1095 \skewchar\elevensy='60
\font@\elevenex=cmex10 scaled 1095
\font@\elevenmsa=msam10 scaled 1095
\font@\elevenmsb=msbm10 scaled 1095
\font@\eleveneufm=eufm10 scaled 1095

\def\elevenpoint{\p@@1.1\p@ \normallineskiplimit\p@@
 \mathsurround\m@ths@r \normalbaselineskip12\p@@
 \abovedisplayskip12\p@@ plus3\p@@ minus9\p@@
 \belowdisplayskip\abovedisplayskip
 \abovedisplayshortskip\z@ plus3\p@@
 \belowdisplayshortskip7\p@@ plus3\p@@ minus4\p@@
 \textonlyfont@\rm\elevenrm \textonlyfont@\it\elevenit
 \textonlyfont@\sl\elevensl \textonlyfont@\bf\elevenbf
 \textonlyfont@\smc\elevensmc \textonlyfont@\tt\eleventt
 \ifsyntax@ \def\big##1{{\hbox{$\left##1\right.$}}}%
  \let\Big\big \let\bigg\big \let\Bigg\big
 \else
  \textfont\z@\elevenrm \scriptfont\z@\eightrm \scriptscriptfont\z@\sixrm
  \textfont\@ne\eleveni \scriptfont\@ne\eighti \scriptscriptfont\@ne\sixi
  \textfont\tw@\elevensy \scriptfont\tw@\eightsy \scriptscriptfont\tw@\sixsy
  \textfont\thr@@\elevenex \scriptfont\thr@@\eightex
	\scriptscriptfont\thr@@\sevenex
  \textfont\itfam\elevenit \scriptfont\itfam\eightit
	\scriptscriptfont\itfam\sevenit
  \textfont\bffam\elevenbf \scriptfont\bffam\eightbf
	\scriptscriptfont\bffam\sixbf
  \textfont\msafam\elevenmsa \scriptfont\msafam\eightmsa
	\scriptscriptfont\msafam\sixmsa
  \textfont\msbfam\elevenmsb \scriptfont\msbfam\eightmsb
	\scriptscriptfont\msbfam\sixmsb
  \textfont\eufmfam\eleveneufm \scriptfont\eufmfam\eighteufm
	\scriptscriptfont\eufmfam\sixeufm
  \setbox\strutbox\hbox{\vrule height8.5\p@@ depth3.5\p@@ width\z@}%
  \setbox\strutbox@\hbox{\lower.5\normallineskiplimit\vbox{%
	\kern-\normallineskiplimit\copy\strutbox}}%
  \setbox\z@\vbox{\hbox{$($}\kern\z@}\b@gsize1.2\ht\z@
  \fi
  \normalbaselines\rm\dotsspace@1.5mu\ex@.2326ex\jot3\ex@}

\def\m@R@f@[#1]{\mathsurzero{
 \s@ct{}{#1}}\wr@@c{\string\Refcd{#1}{\the\pageno}}\B@gr@
 \frenchspacing\rcount\z@\refkey{\k@yf@nt[##1]}\refno{\k@yf@nt[##1]}%
 \widest{AZ}\keyright\let\Key\key\let\refin\relax}
\def\widest#1{\s@twd@\r@f@nd{\r@fk@y{\k@yf@nt#1}\enspace}}
\def\widestno#1{\s@twd@\r@f@nd{\r@fn@{\k@yf@nt#1}\enspace}}
\def\widestlabel#1{\s@twd@\r@f@nd{\k@yf@nt#1\enspace}}
\def\refkey{\def\r@fk@y##1} \def\refno{\def\r@fn@##1}
\def\keyright{\def\r@fit@m{\hang\textindent}}
\def\keyflat{\def\r@fit@m##1{\setbox\z@\hbox{##1\enspace}\hang\noindent
 \ifnum\wd\z@<\parindent\indent\hglue-\wd\z@\fi\unhbox\z@}}

\def\R@fb@x{\global\setbox\r@f@b@x} \def\K@yb@x{\global\setbox\k@yb@x}
\def\ref{\par\b@gr@\r@ff@nt\R@fb@x\box\voidb@x\K@yb@x\box\voidb@x
 \@fn@mfalse\@fl@bfalse\b@g@nr@f}
\def\c@nc@t#1{\setbox\z@\lastbox
 \setbox\adjb@x\hbox{\unhbox\adjb@x\unhbox\z@\unskip\unskip\unpenalty#1}}
\def\adjust#1{\relax\ifmmode\penalty-\@M\null\hfil$\clubpenalty\z@
 \widowpenalty\z@\interlinepenalty\z@\offinterlineskip\endgraf
 \setbox\z@\lastbox\unskip\unpenalty\c@nc@t{#1}\nt$\hfil\penalty-\@M
 \else\endgraf\c@nc@t{#1}\nt\fi}
\def\adjustnext#1{\P@nct\hbox{#1}\ignore}
\def\adjustend#1{\def\@djp@{#1}\ignore}
\def\addtoks#1{\global\@ddt@ks{#1}\ignore}
\def\addnext#1{\global\@dd@p@n{#1}\ignore}

\def\cl@s@{\adjust{\@djp@}\endgraf\setbox\z@\lastbox
 \global\setbox\@ne\hbox{\unhbox\adjb@x\ifvoid\z@\else\unhbox\z@\unskip\unskip
 \unpenalty\fi}\egroup\ifnum\c@rr@nt=\k@yb@x\global\fi
 \setbox\c@rr@nt\hbox{\unhbox\@ne\box\p@nct@}\P@nct\null
 \the\@ddt@ks\global\@ddt@ks{}}
\def\@p@n#1{\def\c@rr@nt{#1}\setbox\c@rr@nt\vbox\bgroup\let\@djp@\relax
 \hsize\maxdimen\nt\the\@dd@p@n\global\@dd@p@n{}}
\def\b@g@nr@f{\bgroup\@p@n\z@}
\def\key{\cl@s@\ifvoid\k@yb@x\@p@n\k@yb@x\k@yf@nt\else\@p@n\z@\fi}
\def\label{\cl@s@\ifvoid\k@yb@x\global\@fl@btrue\@p@n\k@yb@x\k@yf@nt\else
 \@p@n\z@\fi}
\def\no{\cl@s@\ifvoid\k@yb@x\gad\rcount\global\@fn@mtrue
 \K@yb@x\hbox{\k@yf@nt\the\rcount}\fi\@p@n\z@}
\def\labelno{\cl@s@\ifvoid\k@yb@x\gad\rcount\@fl@btrue
 \@p@n\k@yb@x\k@yf@nt\the\rcount\else\@p@n\z@\fi}
\def\by{\cl@s@\@p@n\b@b@x} \def\paper{\cl@s@\@p@n\p@p@rb@x\p@p@rf@nt\ignore}
\def\jour{\cl@s@\@p@n\j@@rb@x} \def\yr{\cl@s@\@p@n\y@@rb@x}
\def\vol{\cl@s@\@p@n\v@lb@x\v@lf@nt\ignore}
\def\issue{\cl@s@\@p@n\is@b@x\iss@f@nt\ignore}
\def\page{\cl@s@\ifp@g@s\@p@n\z@\else\p@g@true\@p@n\p@g@b@x\fi}
\def\pages{\cl@s@\ifp@g@\@p@n\z@\else\p@g@strue\@p@n\p@g@b@x\fi}
\def\inbook{\cl@s@\@p@n\inb@@kb@x}
\def\book{\cl@s@\@p@n\b@@kb@x\b@@kf@nt\ignore}
\def\publ{\cl@s@\@p@n\p@blb@x} \def\publaddr{\cl@s@\@p@n\p@bl@db@x}
\def\ed{\cl@s@\ifed@s\@p@n\z@\else\ed@true\@p@n\ed@b@x\fi}
\def\eds{\cl@s@\ifed@\@p@n\z@\else\ed@strue\@p@n\ed@b@x\fi}
\def\info{\cl@s@\@p@n\inf@b@x} \def\paperinfo{\cl@s@\@p@n\p@p@nf@b@x}
\def\bookinfo{\cl@s@\@p@n\b@@nf@b@x} 
\def\P@nct{\global\setbox\p@nct@} \def\nopunct{\P@nct\box\voidb@x}
\def\p@@@t#1#2{\ifvoid\p@nct@\else#1\unhbox\p@nct@#2\fi}
\def\sp@@{\penalty-50 \space\hskip\z@ plus.1em}
\def\c@mm@{\p@@@t,\sp@@} \def\sp@c@{\p@@@t\empty\sp@@}
\def\p@tb@x#1#2{\ifvoid#1\else#2\@nb@x#1\fi}
\def\@nb@x#1{\unhbox#1\P@nct\lastbox}
\def\endr@f@{\cl@s@\nopunct
 \R@fb@x\hbox{\unhbox\r@f@b@x \p@tb@x\b@b@x\empty
 \ifvoid\j@@rb@x\ifvoid\inb@@kb@x\ifvoid\p@p@rb@x\ifvoid\b@@kb@x
  \ifvoid\p@p@nf@b@x\ifvoid\b@@nf@b@x
  \p@tb@x\v@lb@x\c@mm@ \ifvoid\y@@rb@x\else\sp@c@(\@nb@x\y@@rb@x)\fi
  \p@tb@x\is@b@x\c@mm@ \p@tb@x\p@g@b@x\c@mm@ \p@tb@x\inf@b@x\c@mm@
  \else\p@tb@x \b@@nf@b@x\c@mm@ \p@tb@x\v@lb@x\c@mm@ \p@tb@x\is@b@x\sp@c@
  \ifvoid\ed@b@x\else\sp@c@(\@nb@x\ed@b@x,\space\ifed@ ed.\else eds.\fi)\fi
  \p@tb@x\p@blb@x\c@mm@ \p@tb@x\p@bl@db@x\c@mm@ \p@tb@x\y@@rb@x\c@mm@
  \p@tb@x\p@g@b@x{\c@mm@\ifp@g@ p\p@@nt\else pp\p@@nt\fi}%
  \p@tb@x\inf@b@x\c@mm@\fi
  \else \p@tb@x\p@p@nf@b@x\c@mm@ \p@tb@x\v@lb@x\c@mm@
  \ifvoid\y@@rb@x\else\sp@c@(\@nb@x\y@@rb@x)\fi
  \p@tb@x\is@b@x\c@mm@ \p@tb@x\p@g@b@x\c@mm@ \p@tb@x\inf@b@x\c@mm@\fi
  \else \p@tb@x\b@@kb@x\c@mm@
  \p@tb@x\b@@nf@b@x\c@mm@ \p@tb@x\p@blb@x\c@mm@
  \p@tb@x\p@bl@db@x\c@mm@ \p@tb@x\y@@rb@x\c@mm@
  \ifvoid\p@g@b@x\else\c@mm@\@nb@x\p@g@b@x p\fi \p@tb@x\inf@b@x\c@mm@ \fi
  \else \c@mm@\@nb@x\p@p@rb@x\ic@\p@tb@x\p@p@nf@b@x\c@mm@
  \p@tb@x\v@lb@x\sp@c@ \ifvoid\y@@rb@x\else\sp@c@(\@nb@x\y@@rb@x)\fi
  \p@tb@x\is@b@x\c@mm@ \p@tb@x\p@g@b@x\c@mm@\p@tb@x\inf@b@x\c@mm@\fi
  \else \p@tb@x\p@p@rb@x\c@mm@\ic@\p@tb@x\p@p@nf@b@x\c@mm@
  \c@mm@\@nb@x\inb@@kb@x \p@tb@x\b@@nf@b@x\c@mm@ \p@tb@x\v@lb@x\sp@c@
  \p@tb@x\is@b@x\sp@c@
  \ifvoid\ed@b@x\else\sp@c@(\@nb@x\ed@b@x,\space\ifed@ ed.\else eds.\fi)\fi
  \p@tb@x\p@blb@x\c@mm@ \p@tb@x\p@bl@db@x\c@mm@ \p@tb@x\y@@rb@x\c@mm@
  \p@tb@x\p@g@b@x{\c@mm@\ifp@g@ p\p@@nt\else pp\p@@nt\fi}%
  \p@tb@x\inf@b@x\c@mm@\fi
  \else\p@tb@x\p@p@rb@x\c@mm@\ic@\p@tb@x\p@p@nf@b@x\c@mm@\p@tb@x\j@@rb@x\c@mm@
  \p@tb@x\v@lb@x\sp@c@ \ifvoid\y@@rb@x\else\sp@c@(\@nb@x\y@@rb@x)\fi
  \p@tb@x\is@b@x\c@mm@ \p@tb@x\p@g@b@x\c@mm@ \p@tb@x\inf@b@x\c@mm@ \fi}}
\def\m@r@f#1#2{\endr@f@\ifvoid\p@nct@\else\R@fb@x\hbox{\unhbox\r@f@b@x
 #1\unhbox\p@nct@\penalty-200\enskip#2}\fi\egroup\b@g@nr@f}
\def\endref{\endr@f@\ifvoid\p@nct@\else\R@fb@x\hbox{\unhbox\r@f@b@x.}\fi
 \parindent\r@f@nd
 \r@fit@m{\ifvoid\k@yb@x\else\if@fn@m\r@fn@{\unhbox\k@yb@x}\else
 \if@fl@b\unhbox\k@yb@x\else\r@fk@y{\unhbox\k@yb@x}\fi\fi\fi}\unhbox\r@f@b@x
 \endgraf\egroup\endgroup}
\def\moreref{\m@r@f;\empty}
\def\transl{\m@r@f;{\unskip\space
 {\sl English translation\ic@}:\penalty-66 \space}}
\def\endRefs{\endgraf\goodbreak\endgroup}

\hyphenation{acad-e-my acad-e-mies af-ter-thought anom-aly anom-alies
an-ti-deriv-a-tive an-tin-o-my an-tin-o-mies apoth-e-o-ses
apoth-e-o-sis ap-pen-dix ar-che-typ-al as-sign-a-ble as-sist-ant-ship
as-ymp-tot-ic asyn-chro-nous at-trib-uted at-trib-ut-able bank-rupt
bank-rupt-cy bi-dif-fer-en-tial blue-print busier busiest
cat-a-stroph-ic cat-a-stroph-i-cally con-gress cross-hatched data-base
de-fin-i-tive de-riv-a-tive dis-trib-ute dri-ver dri-vers eco-nom-ics
econ-o-mist elit-ist equi-vari-ant ex-quis-ite ex-tra-or-di-nary
flow-chart for-mi-da-ble forth-right friv-o-lous ge-o-des-ic
ge-o-det-ic geo-met-ric griev-ance griev-ous griev-ous-ly
hexa-dec-i-mal ho-lo-no-my ho-mo-thetic ideals idio-syn-crasy
in-fin-ite-ly in-fin-i-tes-i-mal ir-rev-o-ca-ble key-stroke
lam-en-ta-ble light-weight mal-a-prop-ism man-u-script mar-gin-al
meta-bol-ic me-tab-o-lism meta-lan-guage me-trop-o-lis
met-ro-pol-i-tan mi-nut-est mol-e-cule mono-chrome mono-pole
mo-nop-oly mono-spline mo-not-o-nous mul-ti-fac-eted mul-ti-plic-able
non-euclid-ean non-iso-mor-phic non-smooth par-a-digm par-a-bol-ic
pa-rab-o-loid pa-ram-e-trize para-mount pen-ta-gon phe-nom-e-non
post-script pre-am-ble pro-ce-dur-al pro-hib-i-tive pro-hib-i-tive-ly
pseu-do-dif-fer-en-tial pseu-do-fi-nite pseu-do-nym qua-drat-ic
quad-ra-ture qua-si-smooth qua-si-sta-tion-ary qua-si-tri-an-gu-lar
quin-tes-sence quin-tes-sen-tial re-arrange-ment rec-tan-gle
ret-ri-bu-tion retro-fit retro-fit-ted right-eous right-eous-ness
ro-bot ro-bot-ics sched-ul-ing se-mes-ter semi-def-i-nite
semi-ho-mo-thet-ic set-up se-vere-ly side-step sov-er-eign spe-cious
spher-oid spher-oid-al star-tling star-tling-ly sta-tis-tics
sto-chas-tic straight-est strange-ness strat-a-gem strong-hold
sum-ma-ble symp-to-matic syn-chro-nous topo-graph-i-cal tra-vers-a-ble
tra-ver-sal tra-ver-sals treach-ery turn-around un-at-tached
un-err-ing-ly white-space wide-spread wing-spread wretch-ed
wretch-ed-ly Brown-ian Eng-lish Euler-ian Feb-ru-ary Gauss-ian
Grothen-dieck Hamil-ton-ian Her-mit-ian Jan-u-ary Japan-ese Kor-te-weg
Le-gendre Lip-schitz Lip-schitz-ian Mar-kov-ian Noe-ther-ian
No-vem-ber Rie-mann-ian Schwarz-schild Sep-tem-ber}

\let\nopagenumber\p@gen@false \let\putpagenumber\p@gen@true

\outer\def\myRefs{\myR@fs} \r@st@re\proclaim
\def\bye{\par\vfill\supereject\cl@selbl\cl@secd\b@e} \r@endd@\b@e
 \let\Key\key \def\endpro{\par\endproclaim}
\let\d@c@\document \def\document{\d@c@\tenpoint}
\hyphenation{ortho-gon-al}

\newtoks\@@tp@t \@@tp@t\output
\output=\@ft@{\let\{\noexpand\the\@@tp@t}
\let\{\relax

\newif\ifVersion \Versiontrue
\def\p@n@l#1{\ifnum#1=\z@\else\penalty#1\relax\fi}

\def\s@ct#1#2{\ifVersion
 \skip@\lastskip\ifdim\skip@<1.5\bls\vskip-\skip@\p@n@l{-200}\vsk.5>%
 \p@n@l{-200}\vsk.5>\p@n@l{-200}\vsk.5>\p@n@l{-200}\vsk-1.5>\else
 \p@n@l{-200}\fi\ifdim\skip@<.9\bls\vsk.9>\else
 \ifdim\skip@<1.5\bls\vskip\skip@\fi\fi
 \vtop{\twelvepoint\raggedright\s@cf@nt\vp1\vsk->\vskip.16ex
 \s@twd@\parindent{#1}%
 \ifdim\parindent>\z@\adv\parindent.5em\fi\hang\textindent{#1}#2\strut}
 \else
 \p@sk@p{-200}{.8\bls}\vtop{\s@cf@nt\s@twd@\parindent{#1}%
 \ifdim\parindent>\z@\adv\parindent.5em\fi\hang\textindent{#1}#2\strut}\fi
 \nointerlineskip\nobreak\vtop{\strut}\nobreak\vskip-.6\bls\nobreak}

\def\s@bs@ct#1#2{\ifVersion
 \skip@\lastskip\ifdim\skip@<1.5\bls\vskip-\skip@\p@n@l{-200}\vsk.5>%
 \p@n@l{-200}\vsk.5>\p@n@l{-200}\vsk.5>\p@n@l{-200}\vsk-1.5>\else
 \p@n@l{-200}\fi\ifdim\skip@<.9\bls\vsk.9>\else
 \ifdim\skip@<1.5\bls\vskip\skip@\fi\fi
 \vtop{\elevenpoint\raggedright\s@bf@nt\vp1\vsk->\vskip.16ex%
 \s@twd@\parindent{#1}\ifdim\parindent>\z@\adv\parindent.5em\fi
 \hang\textindent{#1}#2\strut}
 \else
 \p@sk@p{-200}{.6\bls}\vtop{\s@bf@nt\s@twd@\parindent{#1}%
 \ifdim\parindent>\z@\adv\parindent.5em\fi\hang\textindent{#1}#2\strut}\fi
 \nointerlineskip\nobreak\vtop{\strut}\nobreak\vskip-.8\bls\nobreak}

\def\gadv{\global\adv} \def\gad#1{\gadv#1\@ne} \def\gadneg#1{\gadv#1-\@ne}

\newcount\t@@n \t@@n=10 \newbox\testbox

\newcount\Sno \newcount\Lno \newcount\Fno

\def\pr@cl#1{\r@st@re\pr@c@\pr@c@{#1}\global\let\pr@c@\relax}

\def\l@L#1{\l@bel{#1}L} \def\l@F#1{\l@bel{#1}F} \def\<#1>{\l@b@l{#1}F}

\def\tagg#1{\tag"\rlap{\rm(#1)}\kern.01\p@"}
\def\Tag#1{\tag{\l@F{#1}}} \def\Tagg#1{\tagg{\l@F{#1}}}

\def\xspace{\kern.34em}

\def\Th#1{\pr@cl{Theorem\xspace\l@L{#1}}\ignore}
\def\Lm#1{\pr@cl{Lemma\xspace\l@L{#1}}\ignore}
\def\Cr#1{\pr@cl{Corollary\xspace\l@L{#1}}\ignore}
\def\Df#1{\pr@cl{Definition\xspace\l@L{#1}}\ignore}
\def\Cj#1{\pr@cl{Conjecture\xspace\l@L{#1}}\ignore}
\def\Prop#1{\pr@cl{Proposition\xspace\l@L{#1}}\ignore}
\def\Rem{\demo{\sl Remark}} 
\def\Pf#1.{\demo{Proof #1}} \def\epf{\qed\enddemo}

\def\Proof#1.{\demo{\let\{\relax Proof #1}\def\t@st@{#1}%
 \ifx\t@st@\empty\else\xdef\@@wr##1##2##3##4{##1{##2##3{\the\cdn@}{##4}}}%
 \wr@@c@{\the\cdn@}{Proof #1}\@@wr\wr@@c\string\subcd{\the\pageno}\fi\ignore}

\def\Ap@x{Appendix}
\def\Appendix{\Sno=64 \t@@n\@ne \wr@@c{\string\Appencd}
 \def\sf@rm{\char\the\Sno} \def\sf@rm@{\Ap@x\space\sf@rm} \def\sf@rm@@{\Ap@x}
 \def\s@ct@n##1##2{\s@ct\empty{\setbox\z@\hbox{##1}\ifdim\wd\z@=\z@
 \if##2*\sf@rm@@\else\if##2.\sf@rm@@.\else##2\fi\fi\else
 \if##2*\sf@rm@\else\if##2.\sf@rm@.\else\sf@rm@.\enspace##2\fi\fi\fi}}}
\def\Appcd#1#2#3{\gad\Cdentry\global\cdentry\z@\def\Ap@@{\hglue-\l@ftcd\Ap@x}
 \ifx\@ppl@ne\empty\def\l@@b{\@fwd@@{#1}{\space#1}{}}
 \if*#2\entcd{}{\Ap@@\l@@b}{#3}\else\if.#2\entcd{}{\Ap@@\l@@b.}{#3}\else
 \entcd{}{\Ap@@\l@@b.\enspace#2}{#3}\fi\fi\else
 \def\l@@b{\@fwd@@{#1}{\c@l@b{#1}}{}}\if*#2\entcd{\l@@b}{\Ap@x}{#3}\else
 \if.#2\entcd{\l@@b}{\Ap@x.}{#3}\else\entcd{\l@@b}{#2}{#3}\fi\fi\fi}

\let\s@ct@n\s@ct
\def\s@ct@@[#1]#2{\@ft@\xdef\csname @#1@S@\endcsname{\sf@rm}\wr@@x{}%
 \wr@@x{\string\labeldef{S}\space{\?#1@S?}\space{#1}}%
 {
 \s@ct@n{\sf@rm@}{#2}}\wr@@c{\string\Entcd{\?#1@S?}{#2}{\the\pageno}}}
\def\s@ct@#1{\wr@@x{}{
 \s@ct@n{\sf@rm@}{#1}}\wr@@c{\string\Entcd{\sf@rm}{#1}{\the\pageno}}}
\def\s@ct@e[#1]#2{\@ft@\xdef\csname @#1@S@\endcsname{\sf@rm}\wr@@x{}%
 \wr@@x{\string\labeldef{S}\space{\?#1@S?}\space{#1}}%
 {
 \s@ct@n\empty{#2}}\wr@@c{\string\Entcd{}{#2}{\the\pageno}}}
\def\s@cte#1{\wr@@x{}{
 \s@ct@n\empty{#1}}\wr@@c{\string\Entcd{}{#1}{\the\pageno}}}
\def\theSno#1#2{\dff\?#1@S?{#2}%
 \wr@@x{\string\labeldef{S}\space{#2}\space{#1}}\fi}

\newif\ifd@bn@\d@bn@true
\def\Section{\gad\Sno\ifd@bn@\Fno\z@\Lno\z@\fi\@fn@xt[\s@ct@@\s@ct@}
\def\section{\gad\Sno\ifd@bn@\Fno\z@\Lno\z@\fi\@fn@xt[\s@ct@e\s@cte}
 
\def\subsection{\@fn@xt*\subs@ct@\subs@ct}
\def\subs@ct#1{{\s@bs@ct\empty{#1}}\wr@@c{\string\subcd{#1}{\the\pageno}}}
\def\subs@ct@*#1{\vsk->\nobreak
 {\s@bs@ct\empty{#1}}\wr@@c{\string\subcd{#1}{\the\pageno}}}

\def\l@b@l#1#2{\def\n@@{\csname #2no\endcsname}%
 \if*#1\gad\n@@ \@ft@\xdef\csname @#1@#2@\endcsname{\l@f@rm}\else\def\t@st{#1}%
 \ifx\t@st\empty\gad\n@@ \@ft@\xdef\csname @#1@#2@\endcsname{\l@f@rm}%
 \else\@ft@\ifx\csname @#1@#2@mark\endcsname\relax\gad\n@@
 \@ft@\xdef\csname @#1@#2@\endcsname{\l@f@rm}%
 \@ft@\gdef\csname @#1@#2@mark\endcsname{}%
 \wr@@x{\string\labeldef{#2}\space{\?#1@#2?}\space\ifnum\n@@<10 \space\fi{#1}}%
 \fi\fi\fi}
\def\labeldef#1#2#3{\dff\?#3@#1?{#2}}
\def\Labeldef#1#2#3{\dff\?#3@#1?{#2}\@ft@\gdef\csname @#3@#1@mark\endcsname{}}

\def\l@bel#1#2{\l@b@l{#1}{#2}\?#1@#2?}

\newcount\c@cite
\def\?#1?{\csname @#1@\endcsname}
\def\[{\@fn@xt:\c@t@sect\c@t@}
\def\c@t@#1]{{\c@cite\z@\@fwd@@{\?#1@L?}{\adv\c@cite1}{}%
 \@fwd@@{\?#1@F?}{\adv\c@cite1}{}\@fwd@@{\?#1?}{\adv\c@cite1}{}%
 \relax\ifnum\c@cite=\z@{\bf ???}\wrs@x{No label [#1]}\else
 \ifnum\c@cite=1\let\@@PS\relax\let\@@@\relax\else\let\@@PS\underbar
 \def\@@@{{\rm<}}\fi\@@PS{\?#1?\@@@\?#1@L?\@@@\?#1@F?}\fi}}
\def\(#1){{\rm(\c@t@#1])}}
\def\c@t@s@ct#1{\@fwd@@{\?#1@S?}{\?#1@S?\relax}%
 {{\bf ???}\wrs@x{No section label {#1}}}}
\def\c@t@sect:#1]{\c@t@s@ct{#1}} \let\SNo\c@t@s@ct

\newdimen\l@ftcd \newdimen\r@ghtcd \let\nlc\relax
\newcount\Cdentry \newcount\cdentry \let\prentcd\relax \let\postentcd\relax

\def\d@tt@d{\leaders\hbox to 1em{\kern.1em.\hfil}\hfill}
\def\entcd#1#2#3{\gad\cdentry\prentcd\item{\l@bcdf@nt#1}{\entcdf@nt#2}\alb
 \kern.9em\hbox{}\kern-.9em\d@tt@d\kern-.36em{\p@g@cdf@nt#3}\kern-\r@ghtcd
 \hbox{}\postentcd\par}
\def\Entcd#1#2#3{\gad\Cdentry\global\cdentry\z@
 \def\l@@b{\@fwd@@{#1}{\c@l@b{#1}}{}}\vsk.2>\entcd{\l@@b}{#2}{#3}}
\def\subcd#1#2{{\adv\leftskip.333em\entcd{}{\s@bcdf@nt#1}{#2}}}
\def\Refcd#1#2{\def\t@@st{#1}\ifx\t@@st\empty\ifx\r@fl@ne\empty\relax\else
 \R@fcd{\r@fl@ne}{#2}\fi\else\R@fcd{#1}{#2}\fi}
\def\R@fcd#1#2{\sk@@p{.6\bls}\entcd{}{\hglue-\l@ftcd\R@fcdf@nt #1}{#2}}
\def\Refline{\def\r@fl@ne} \def\Refempty{\let\r@fl@ne\empty}
\def\Appencd{\par\adv\leftskip-\l@ftcd\adv\rightskip-\r@ghtcd\@ppl@ne
 \adv\leftskip\l@ftcd\adv\rightskip\r@ghtcd\let\Entcd\Appcd}
\def\appline{\def\@ppl@ne} \def\Appempty{\let\@ppl@ne\empty}
\def\Appline#1{\def\@ppl@ne{\s@bs@ct{}{#1}}}
\def\Leftcd#1{\adv\leftskip-\l@ftcd\s@twd@\l@ftcd{\c@l@b{#1}\enspace}
 \adv\leftskip\l@ftcd}
\def\Rightcd#1{\adv\rightskip-\r@ghtcd\s@twd@\r@ghtcd{#1\enspace}
 \adv\rightskip\r@ghtcd}
\def\C@nt{Contents} \def\Ap@s{Appendices} \def\R@fcs{References}
\def\contents{\@fn@xt*\cont@@\cont@}
\def\cont@{\@fn@xt[\cnt@{\cnt@[\C@nt]}}
\def\cont@@*{\@fn@xt[\cnt@@{\cnt@@[\C@nt]}}
\def\cnt@[#1]{\c@nt@{M}{#1}{44}{\s@bs@ct{}{\@ppl@f@nt\Ap@s}}}
\def\cnt@@[#1]{\c@nt@{M}{#1}{44}{}}
\def\endco{\par\penalty-500\vsk>\vskip\z@\endgroup}
\def\readcd{\@np@t{\jobname.cd}}
\def\Cde{\@fn@xt*\Cde@@\Cde@}
\def\Cde@{\@fn@xt[\Cd@{\Cd@[\C@nt]}}
\def\Cde@@*{\@fn@xt[\Cd@@{\Cd@@[\C@nt]}}
\def\Cd@[#1]{\cnt@[#1]\readcd\endco}
\def\Cd@@[#1]{\cnt@@[#1]\readcd\endco}
\def\contlabeldef{\def\c@l@b}

\long\def\c@nt@#1#2#3#4{\s@twd@\l@ftcd{\c@l@b{#1}\enspace}
 \s@twd@\r@ghtcd{#3\enspace}\adv\r@ghtcd1.333em
 \def\@ppl@ne{#4}\def\r@fl@ne{\R@fcs}\s@ct{}{#2}\B@gr@\parindent\z@\let\nlc\nl
 \let\nl\relax\parskip.2\bls\adv\leftskip\l@ftcd\adv\rightskip\r@ghtcd}

\def\writecd{\immediate\openout\@@cd\jobname.cd \def\wr@@c{\write\@@cd}
 \def\cl@secd{\immediate\write\@@cd{\string\endinput}\immediate\closeout\@@cd}
 \def\closecd{\cl@secd\global\let\cl@secd\relax}}
\let\cl@secd\relax \def\wr@@c#1{} \let\closecd\relax

\def\dff{\@ft@\d@f} \def\d@f{\@ft@\def}
\def\edff{\@ft@\ed@f} \def\ed@f{\@ft@\edef}
\def\gdff{\@ft@\gd@f} \def\gd@f{\@ft@\gdef}
\def\defi#1#2{\def#1{#2}\wr@@x{\string\def\string#1{#2}}}

\def\qed{\hbox{}\nobreak\hfill\nobreak{\m@th$\,\square$}}
\def\back#1 {\strut\kern-.33em #1\enspace\ignore} 

\def\hcor#1{\advance\hoffset by #1}
\def\vcor#1{\advance\voffset by #1}
\let\bls\baselineskip \let\ignore\ignorespaces
\ifx\ic@\undefined \let\ic@\/\fi
\def\vsk#1>{\vskip#1\bls} \let\adv\advance
\def\vv#1>{\vadjust{\vsk#1>}\ignore}
\def\vvn#1>{\vadjust{\nobreak\vsk#1>\nobreak}\ignore}
\def\vvv#1>{\vskip\z@\vsk#1>\nt\ignore}
\def\vvgood{\vadjust{\penalty-500}}
\def\nngood{\noalign{\penalty-500}}

\def\Goodbreak{\par\penalty-\@m}
\def\wgood#1>{\vv#1>\vvgood\vv-#1>}
\def\wwgood#1:#2>{\vv#1>\vvgood\vv#2>}
\def\mmgood#1:#2>{\cnn#1>\nngood\cnn#2>}
\def\goodsk#1:#2>{\vsk#1>\goodbreak\vsk#2>\vsk0>}
\def\ragood{\vadjust{\vskip\z@ plus 12pt}\vvgood}

\def\Par{\vsk.5>} \def\setparindent{\edef\Parindent{\the\parindent}}
\def\Type{\vsk.5>\bgroup\parindent\z@\tt\rightskip\z@ plus1em minus1em%
 \spaceskip.3333em \xspaceskip.5em\relax}
\def\endType{\vsk.5>\egroup\nt} 

\let\Hat\widehat \let\Tilde\widetilde \let\dollar\$ \let\ampersand\&
\let\sss\scriptscriptstyle  
\let\vp\vphantom \let\hp\hphantom \let\nt\noindent
\let\cline\centerline \let\lline\leftline \let\rline\rightline
\def\nn#1>{\noalign{\vskip#1\p@@}} \def\NN#1>{\openup#1\p@@}
\def\cnn#1>{\noalign{\vsk#1>}}
 
\let\Lim\lim \def\lim{\Lim\limits} \let\Sum\sum \def\sum{\Sum\limits}
\def\Plus{\bigoplus\limits} 
\let\Prod\prod \def\prod{\Prod\limits} \let\Int\int \def\int{\Int\limits}

\def\tsum{\mathop{\tsize\Sum}\limits} 
 \def\&{.\kern.1em}
\def\nl{\leavevmode\hfill\break} \def\~{\leavevmode\@fn@xt~\m@n@s\@md@@sh}
\def\@md@@sh{\@fn@xt-\d@@sh\@md@sh} \def\@md@sh{\raise.13ex\hbox{--}}
\def\m@n@s~{\raise.15ex\mbox{-}} \def\d@@sh-{\raise.15ex\hbox{-}}

\let\procent\% \def\%#1{\ifmmode\mathop{#1}\limits\else\procent#1\fi}
\let\@ml@t\" \def\"#1{\ifmmode ^{(#1)}\else\@ml@t#1\fi}
\let\@c@t@\' \def\'#1{\ifmmode _{(#1)}\else\@c@t@#1\fi}
\let\colon\: \def\:{^{\vp{\topsmash|}}} 

\let\texspace\ \def\ {\ifmmode\alb\fi\texspace} \def\.{\d@t\ignore}

\newif\ifNewskips

\def\Newskips{\global\Newskipstrue
 \gdef\>{\RIfM@\mskip.666667\thinmuskip\relax\else\kern.111111em\fi}
 \gdef\}{\RIfM@\mskip-.666667\thinmuskip\relax\else\kern-.111111em\fi}
 \gdef\){\RIfM@\mskip.333333\thinmuskip\relax\else\kern.0555556em\fi}
 \gdef\]{\RIfM@\mskip-.333333\thinmuskip\relax\else\kern-.0555556em\fi}}
\def\d@t{\ifNewskips.\hskip.3em\else\def\d@t{.\ }\fi} \def\.{\d@t\ignore}
\Newskips

\let\n@wp@ge\newpage \def\newpage{\endgraf\n@wp@ge}
\let\=\m@th \def\mbox#1{\hbox{\m@th$#1$}}
\def\mtext#1{\text{\m@th$#1$}} \def\^#1{\text{\m@th#1}}
\def\Line#1{\kern-.5\hsize\line{\m@th$\dsize#1$}\kern-.5\hsize}
\def\Lline#1{\kern-.5\hsize\lline{\m@th$\dsize#1$}\kern-.5\hsize}
\def\Cline#1{\kern-.5\hsize\cline{\m@th$\dsize#1$}\kern-.5\hsize}
\def\Rline#1{\kern-.5\hsize\rline{\m@th$\dsize#1$}\kern-.5\hsize}

\def\Ll@p#1{\llap{\m@th$#1$}} \def\Rl@p#1{\rlap{\m@th$#1$}}
 \def\Cl@p#1{\llap{\m@th$#1$\hss}}
\def\Llap#1{\mathchoice{\Ll@p{\dsize#1}}{\Ll@p{\tsize#1}}{\Ll@p{\ssize#1}}%
 {\Ll@p{\sss#1}}}
\def\Clap#1{\mathchoice{\Cl@p{\dsize#1}}{\Cl@p{\tsize#1}}{\Cl@p{\ssize#1}}%
 {\Cl@p{\sss#1}}}
\def\Rlap#1{\mathchoice{\Rl@p{\dsize#1}}{\Rl@p{\tsize#1}}{\Rl@p{\ssize#1}}%
 {\Rl@p{\sss#1}}}
 
\def\LRtph#1#2{\setbox\z@\hbox{#1}\dimen\z@\wd\z@\hbox{\hbox to\dimen\z@{#2}}}
\def\LRph#1#2{\LRtph{\m@th$#1$}{\m@th$#2$}}

\def\LLph#1#2{\LRph{#1}{\hss#2}} 

\def\Lph#1#2{\mathchoice{\LLph{\dsize#1}{\dsize#2}}{\LLph{\tsize#1}{\tsize#2}}
 {\LLph{\ssize#1}{\ssize#2}}{\LLph{\sss#1}{\sss#2}}}

\def\Lto#1{\setbox\z@\mbox{\tsize{#1}}%
 \mathrel{\mathop{\hbox to\wd\z@{\rightarrowfill}}\limits#1}}
\def\Lgets#1{\setbox\z@\mbox{\tsize{#1}}%
 \mathrel{\mathop{\hbox to\wd\z@{\leftarrowfill}}\limits#1}}
 \def\vpp#1{{\vp{\big]}}_{#1}}

\let\alb\allowbreak

 \let\x\times \let\ox\otimes 
  \let\tabs\+
 \let\ge\geqslant
\let\der\partial \let\8\infty \let\*\star
 \let\ket\rangle

 \def\vert{\ |\ }

\let\lb\lbrace \let\rb\rbrace

\def\lsym#1{#1\alb\ldots\relax#1\alb}
\def\lc{\lsym,}   \def\lox{\lsym\ox}

\def\End{\mathop{\roman{End}\>}\nolimits}

 \def\1{^{-1}} \let\underscore\_ \def\_#1{_{\Rlap{#1}}}
\def\vst#1{{\lower1.9\p@@\mbox{\bigr|_{\raise.5\p@@\mbox{\ssize#1}}}}}
\def\vrp#1:#2>{{\vrule height#1 depth#2 width\z@}}
\def\vru#1>{\vrp#1:\z@>} \def\vrd#1>{\vrp\z@:#1>}
\def\qqq{\qquad\quad} 
\def\sscr#1{\raise.3ex\mbox{\sss#1}} \def\@@PS{\bold{OOPS!!!}}

\def\intcl{\mathop
 {\Rlap{\raise.3ex\mbox{\kern.12em\curvearrowleft}}\int}\limits}
\def\intcr{\mathop
 {\Rlap{\raise.3ex\mbox{\kern.24em\curvearrowright}}\int}\limits}

\def\pms{\raise.25ex\mbox{\ssize\pm}\>}
\def\mps{\raise.25ex\mbox{\ssize\mp}\>}

\let\gm\gamma  
  
 \let\eps\varepsilon \let\epsilon\eps

\let\ka\kappa

 \let\Si\Sigma
 
 \let\phi\varphi

\def\C{\Bbb C}

\def\Z{\Bbb Z}

\def\BB{\Bbb B}

 \def\Zpp{\Z_{>0}}

\def\difl/{differential} \def\dif/{difference}
\def\cf.{cf.\ \ignore} \def\Cf.{Cf.\ \ignore}
\def\egv/{eigenvector} \def\eva/{eigenvalue} \def\eq/{equation}
\def\lhs/{the left hand side} \def\rhs/{the right hand side}
\def\Lhs/{The left hand side} \def\Rhs/{The right hand side}
\def\gby/{generated by} \def\wrt/{with respect to} \def\st/{such that}
\def\resp/{respectively} \def\off/{offdiagonal} \def\wt/{weight}
\def\pol/{polynomial} \def\rat/{rational} \def\tri/{trigonometric}
\def\fn/{function} \def\var/{variable} \def\raf/{\rat/ \fn/}
\def\inv/{invariant} \def\hol/{holomorphic} \def\hof/{\hol/ \fn/}
\def\mer/{meromorphic} \def\mef/{\mer/ \fn/} \def\mult/{multiplicity}
\def\sym/{symmetric} \def\perm/{permutation} \def\fd/{finite-dimensional}
\def\rep/{representation} \def\irr/{irreducible} \def\irrep/{\irr/ \rep/}
\def\hom/{homomorphism} \def\aut/{automorphism} \def\iso/{isomorphism}
\def\lex/{lexicographical} \def\as/{asymptotic} \def\asex/{\as/ expansion}
\def\ndeg/{nondegenerate} \def\neib/{neighbourhood} \def\deq/{\dif/ \eq/}
\def\hw/{highest \wt/} \def\gv/{generating vector} \def\eqv/{equivalent}
\def\msd/{method of steepest descend} \def\pd/{pairwise distinct}
\def\wlg/{without loss of generality} \def\Wlg/{Without loss of generality}
\def\onedim/{one-dimensional} \def\qcl/{quasiclassical} \def\hwv/{\hw/ vector}
\def\hgeom/{hyper\-geometric} \def\hint/{\hgeom/ integral}
\def\hwm/{\hw/ module} \def\emod/{evaluation module} \def\Vmod/{Verma module}
\def\symg/{\sym/ group} \def\sol/{solution} \def\eval/{evaluation}
\def\anf/{analytic \fn/} \def\anco/{analytic continuation}
\def\qg/{quantum group} \def\qaff/{quantum affine algebra}

\def\Rm/{\^{$R$-}matrix} \def\Rms/{\^{$R$-}matrices} \def\YB/{Yang-Baxter \eq/}
\def\Ba/{Bethe ansatz} \def\Bv/{Bethe vector} \def\Bae/{\Ba/ \eq/}
\def\KZv/{Knizh\-nik-Zamo\-lod\-chi\-kov} \def\KZvB/{\KZv/-Bernard}
\def\KZ/{{\sl KZ\/}} \def\qKZ/{{\sl qKZ\/}}
\def\KZB/{{\sl KZB\/}} \def\qKZB/{{\sl qKZB\/}}
\def\qKZo/{\qKZ/ operator} \def\qKZc/{\qKZ/ connection}
\def\KZe/{\KZ/ \eq/} \def\qKZe/{\qKZ/ \eq/} \def\qKZBe/{\qKZB/ \eq/}
\def\XXX/{{\sl XXX\/}} \def\XXZ/{{\sl XXZ\/}} \def\XYZ/{{\sl XYZ\/}}

\def\h@ph{\discretionary{}{}{-}} \def\$#1$-{\,\^{$#1$}\h@ph}

\def\TFT/{Research Insitute for Theoretical Physics}
\def\HY/{University of Helsinki} \def\AoF/{the Academy of Finland}
\def\CNRS/{Supported in part by MAE\~MICECO\~CNRS Fellowship}
\def\LPT/{Laboratoire de Physique Th\'eorique ENSLAPP}
\def\ENSLyon/{\'Ecole Normale Sup\'erieure de Lyon}
\def\LPTaddr/{46, All\'ee d'Italie, 69364 Lyon Cedex 07, France}
\def\enslapp/{URA 14\~36 du CNRS, associ\'ee \`a l'E.N.S.\ de Lyon,
au LAPP d'Annecy et \`a l'Universit\`e de Savoie}
\def\ensemail/{vtarasov\@ enslapp.ens-lyon.fr}
\def\DMS/{Department of Mathematics, Faculty of Science}
\def\DMO/{\DMS/, Osaka University}
\def\DMOaddr/{Toyonaka, Osaka 560, Japan}
\def\dmoemail/{vt\@ math.sci.osaka-u.ac.jp}
\def\MPI/{Max\)-Planck\)-Institut} \def\MPIM/{\MPI/ f\"ur Mathematik}
\def\MPIMaddr/{P\]\&O.\ Box 7280, D\~-\]53072 \,Bonn, Germany}
\def\mpimemail/{tarasov\@ mpim-bonn.mpg.de}
\def\SPb/{St\&Peters\-burg}
\def\home/{\SPb/ Branch of Steklov Mathematical Institute}
\def\homeaddr/{Fontanka 27, \SPb/ \,191011, Russia}
\def\homemail/{vt\@ pdmi.ras.ru}
\def\absence/{On leave of absence from \home/}
\def\support/{Supported in part by}
\def\UNC/{Department of Mathematics, University of North Carolina}
\def\ChH/{Chapel Hill} \def\UNCaddr/{\ChH/, NC 27599, USA}
\def\avemail/{anv\@ email.unc.edu}	
\def\grant/{NSF grant DMS\~9501290}	
\def\Grant/{\support/ \grant/}

\def\Aomoto/{K\&Aomoto}
\def\Cher/{I\&Che\-red\-nik}
\def\Dri/{V\]\&G\&Drin\-feld}
\def\Fadd/{L\&D\&Fad\-deev}
\def\Feld/{G\&Felder}
\def\Fre/{I\&B\&Fren\-kel}
\def\Etingof/{P\]\&Etingof}
\def\Gustaf/{R\&A\&Gustafson}
\def\Izergin/{A\&G\&Izergin}
\def\Jimbo/{M\&Jimbo}
\def\Kazh/{D\&Kazhdan}
\def\Kor/{V\]\&E\&Kore\-pin}
\def\Kulish/{P\]\&P\]\&Ku\-lish}
\def\Lusz/{G\&Lusztig}
\def\Miwa/{T\]\&Miwa}
\def\MN/{M\&Naza\-rov}
\def\Reshet/{N\&Reshe\-ti\-khin} \def\Reshy/{N\&\]Yu\&Reshe\-ti\-khin}
\def\SchV/{V\]\&\]V\]\&Schecht\-man} \def\Sch/{V\]\&Schecht\-man}
\def\Skl/{E\&K\&Sklya\-nin}
\def\Smirnov/{F\]\&A\&Smir\-nov}
\def\Takh/{L\&A\&Takh\-tajan}
\def\VT/{V\]\&Ta\-ra\-sov} \def\VoT/{V\]\&O\&Ta\-ra\-sov}
\def\Varch/{A\&\]Var\-chenko} \def\Varn/{A\&N\&\]Var\-chenko}
\def\Zhel/{D\&P\]\&Zhe\-lo\-ben\-ko}

\def\AiA/{Al\-geb\-ra i Ana\-liz}
\def\DAN/{Do\-kla\-dy AN SSSR}
\def\FAA/{Funk\.Ana\-liz i ego pril.}
\def\Izv/{Iz\-ves\-tiya AN SSSR, ser\.Ma\-tem.}
\def\TMF/{Teo\-ret\.Ma\-tem\.Fi\-zi\-ka}
\def\UMN/{Uspehi Matem.\ Nauk}

\def\AMS/{Amer\.Math\.Society}
\def\AMSa/{AMS \publaddr Providence RI}
\def\AMST/{\AMS/ Transl.,\ Ser\&\)2}
\def\AMSTr/{\AMS/ Transl.,} \def\Ser2{Ser\&\)2}
\def\Astq/{Ast\'erisque}
\def\ContM/{Contemp\.Math.}
\def\CMP/{Comm\.Math\.Phys.}
\def\DMJ/{Duke\.Math\.J.}
\def\Inv/{Invent\.Math.} 
\def\IMRN/{Int\.Math\.Res.\ Notices}
\def\JMP/{J\.Math\.Phys.}
\def\JPA/{J\.Phys.\ A}
\def\JSM/{J\.Soviet Math.}
\def\LMJ/{Leningrad Math.\ J.}
\def\LpMJ/{\SPb/ Math.\ J.}
\def\LMP/{Lett\.Math\.Phys.}
\def\NMJ/{Nagoya Math\.J.}
\def\Nucl/{Nucl\.Phys.\ B}
\def\OJM/{Osaka J\.Math.}
\def\RIMS/{Publ\.RIMS, Kyoto Univ.}
\def\SIAM/{SIAM J\.Math\.Anal.}
\def\SMNS/{Selecta Math., New Series}
\def\TMP/{Theor\.Math\.Phys.}
\def\ZNS/{Zap\. nauch\. semin. LOMI}

\def\ASMP/{Advanced Series in Math.\ Phys.{}}

\def\Birk/{Birkh\"auser}
\def\CUP/{Cambridge University Press} \def\CUPa/{\CUP/ \publaddr Cambridge}
\def\Spri/{Springer\)-Verlag} \def\Spria/{\Spri/ \publaddr Berlin}
\def\WS/{World Scientific} \def\WSa/{\WS/ \publaddr Singapore}

\newbox\lefthbox \newbox\righthbox

\let\sectsep. \let\labelsep. \let\contsep. \let\labelspace\relax
\let\sectpre\relax \let\contpre\relax
\def\sf@rm{\the\Sno} \def\sf@rm@{\sectpre\sf@rm\sectsep}
\def\c@l@b#1{\contpre#1\contsep}
\def\l@f@rm{\ifd@bn@\sf@rm\labelsep\fi\labelspace\the\n@@}

\def\sectformdef{\def\sf@rm}

\let\DoubleNum\d@bn@true \let\SingleNum\d@bn@false

\def\NoNewNum{\let\writeldf\relax\def\l@b@l##1##2{\if*##1%
 \@ft@\xdef\csname @##1@##2@\endcsname{\mbox{*{*}*}}\fi}}
\def\NoNewTime{\def\todaydef##1{\def\today{##1}}
 \def\nowtimedef##1{\def\nowtime{##1}}}
\def\NoInput{\let\Input\input\let\writeldf\relax}
\def\Fixed{\NoNewTime\NoInput}

\newbox\dtlb@x
\def\DateTimeLabel{\global\setbox\dtlb@x\vbox to\z@{\ifMag\eightpoint\else
 \ninepoint\fi\sl\vss\rline\today\rline\nowtime}
 \global\headline{\hfil\box\dtlb@x}}

\def\sectfont#1{\def\s@cf@nt{#1}} \sectfont\bf
\def\subsectfont#1{\def\s@bf@nt{#1}} \subsectfont\it
\def\Entcdfont#1{\def\entcdf@nt{#1}} \Entcdfont\relax
\def\labelcdfont#1{\def\l@bcdf@nt{#1}} \labelcdfont\relax
\def\pagecdfont#1{\def\p@g@cdf@nt{#1}} \pagecdfont\relax
\def\subcdfont#1{\def\s@bcdf@nt{#1}} \subcdfont\it
\def\applefont#1{\def\@ppl@f@nt{#1}} \applefont\bf
\def\Refcdfont#1{\def\R@fcdf@nt{#1}} \Refcdfont\bf

\def\reffont#1{\def\r@ff@nt{#1}} \reffont\rm
\def\keyfont#1{\def\k@yf@nt{#1}} \keyfont\rm
\def\paperfont#1{\def\p@p@rf@nt{#1}} \paperfont\it
\def\bookfont#1{\def\b@@kf@nt{#1}} \bookfont\it
\def\volfont#1{\def\v@lf@nt{#1}} \volfont\bf
\def\issuefont#1{\def\iss@f@nt{#1}} \issuefont{no\p@@nt}

\def\adjustmid#1{\kern-#1\p@\alb\hskip#1\p@\relax}
\def\adjustend#1{\adjustnext{\kern-#1\p@\alb\hskip#1\p@}}

\newif\ifcd 

\tenpoint

\csname beta.def\endcsname
\Fixed

\Magset
\PaperA4

\overfullrule 2pt

\def\Shead{\gad\Sno\Fno0\Lno0\vsk.3>{\elevenpoint\bf\the\Sno.\enspace}\ignore}

\def\Ac{\Cal A}
\def\Dc{\Cal D}

\def\Hs{H_{\sssize\Si}}

\def\q{q^{-1}}
\def\qH{q^{\)H}}
\def\qmH{q^{-H}}
\def\qq{q^{\)2}}
\def\qo{q_{\)0}\:}
\def\qz{q_{\)0}}
\def\qqz{q_{\)0}^{\)2}}

\def\qnum#1{[\)#1\)]\vpp q}
\def\cyc#1{#1,\)#1\qq\lc\)#1q^{\)2M-\)2}}

\def\bv#1{|\)#1\)\ket}
\def\bbv#1{|\}|\)#1\)\ket\!\]\ket}
\def\bbbv#1{|\}|\}|\)#1\)\ket\!\]\ket\!\]\ket}
\def\vox{v_0\]\lox v_0}
\def\Vox{V_{\ell_1}\}\lox V_{\ell_N}}

\def\cm{c_1\lc c_m}
\def\elln{\ell_1\lc\ell_N}
\def\ellsum{\ell_1\]\lsym+\ell_N}

\def\tk{t_1\lc t_k}
\def\tka{t_1\lc\!\Hat{\,t_a}\)\lc t_k}
\def\umi{u_1\lc\!\Hat{\>u_i}\)\lc u_m}
\def\um{u_1\lc u_m}
\def\xm{x_1\lc x_M}
\def\xim{x_{i1}\lc x_{iM}}
\def\Xm{X_1\lc X_m}
\def\Xmi{X_1\lc\!\Hat{\>X_i}\)\lc X_m}
\def\zn{z_1\lc z_N}
\def\zlpn{q^{\)2\ell_1}\]z_1\lc q^{\)2\ell_N}\}z_N}
\def\zlmn{q^{-2\ell_1}\]z_1\lc q^{-2\ell_N}\}z_N}
\def\zlpmn{q^{\)\pm\)2\ell_1}\]z_1\lc q^{\)\pm\)2\ell_N}\}z_N}

\def\Uq{U_q(\frak{sl}_2)} \def\slt{\}\Tilde{\>\frak{sl}_2}}

\def\Lop/{\$L\)$-operator}

\hyphenation{an-iso-tro-py}

\labeldef{F} {1\labelsep \labelspace 1}  {Vl}
\labeldef{F} {1\labelsep \labelspace 2}  {B}
\labeldef{F} {1\labelsep \labelspace 3}  {AB}
\labeldef{F} {1\labelsep \labelspace 4}  {DB}
\labeldef{F} {1\labelsep \labelspace 5}  {AD}
\labeldef{F} {1\labelsep \labelspace 6}  {DA}
\labeldef{F} {1\labelsep \labelspace 7}  {A}
\labeldef{F} {1\labelsep \labelspace 8}  {D}
\labeldef{F} {1\labelsep \labelspace 9}  {xxz}
\labeldef{L} {1\labelsep \labelspace 1}  {BA}
\labeldef{F} {1\labelsep \labelspace 10} {Tbv}
\labeldef{L} {1\labelsep \labelspace 2}  {2}
\labeldef{L} {1\labelsep \labelspace 3}  {3}
\labeldef{L} {1\labelsep \labelspace 4}  {TV}
\labeldef{F} {2\labelsep \labelspace 1}  {Tbbv}
\labeldef{F} {3\labelsep \labelspace 1}  {limB}
\labeldef{F} {3\labelsep \labelspace 2}  {xir}
\labeldef{L} {3\labelsep \labelspace 1}  {BAM}
\labeldef{F} {3\labelsep \labelspace 3}  {A2}
\labeldef{F} {3\labelsep \labelspace 4}  {D2}
\labeldef{F} {3\labelsep \labelspace 5}  {xeq}

\document

\center
{\twelvepoint
\bf
On Bethe vectors for the \XXZ/ model at roots of unity}
\vsk1.2>
{\elevenpoint\smc\VT/}
\vsk>
\it\home/\\\homeaddr/
\endcenter
\ftext{\vsk.24>\nt{\tenpoint\sl E-mail\/{\rm:}\enspace\homemail/}\nl
\ninepoint Partially supported by RFFI grant 02\)\~\>01\~\>00085a}
\vsk2>

Recently in a series of papers \cite{DFM}\), \cite{FM1}\), \cite{FM2}
the \Bae/s and \Bv/s for the six-vertex model with the anisotropy commensurable
with $\pi$ or, as usually said, at roots of unity, were studied. In that case
the spectrum of the transfer-matrix, which is a generating \fn/ of the standard
commuting conservation laws in the model, becomes highly degenerate.
In \cite{FM2} a construction of creation operators, responsible for appearance
of the \Bv/s with the same \eva/s of the transfer-matrix, was suggested in
the framework of the algebraic \Ba/. In the note we extend that construction
to the case of the inhomogeneous arbitrary spin \XXZ/ model. Even for the case
of six-vertex model the proof of the main formulae given in the note is simpler
than the original proof in \cite{FM2}.
\vsk.1>
The detailed exposition of the algebraic \Ba/ method can be found
in \cite{KBI}\).
\vsk.1>
The notation used in the note does not coincide with those of \cite{FM2}
and \cite{KBI}, however a reader can easily establish the correspondence.

\Shead
Consider the inhomogeneous \XXZ/ model on the $N\!$-\)vertex lattice
with the anisotropy $\gm$ and the quasiperiodic boundary conditions.
Let $q=e^{i\gm}\}$. Let $\elln$ be the spins of \rep/s at vertices, $\zn$ \}---
the inhomogeneity parameters, and $\ka$ \}--- the quasiperiodicity parameter,
the periodic boundary conditions corresponding to $\ka=1$. We assume that
\>$\qq\}\ne 1$ and $z_i\ne 0$, $2\ell_i\in\Zpp$ for all $i=1\lc N$.
\vsk.3>
For any $x$ set \,$q^{\)x}\}=\)e^{i\gm x}$ and \,$\dsize
\qnum x\,=\;{q^{\)x}\}-q^{-x}\over q-\q}\;=\;{\sin\)(\gm\)x)\over\sin\gm}\>$.
\vsk.3>
The \XXZ/ model is described by the \Lop/
\vvn.1>
$$
L(u)\,=\>
\pmatrix \>u\)\qH\}-\qmH & u\>(q-\q)\)F\>
\\\nn4>
\>(q-\q)\)E & u\)\qmH\}-\qH\> \endpmatrix \>
\vv.1>
$$
where the elements $E,\,F,\,H$ are generators of $\Uq$:
$$
[\)H\),E\>]\>=\>E\,,\qquad [\)H\),F\>]\>=\>-\)F\>,\qquad
[\)E\),F\>]\>=\>\qnum{2\)H}\,.
\vv-.2>
$$
\vv-.1>
Define a \rep/ of $\Uq$ of nonnegative integral or semiintegral spin $\ell$
in the space \,$V_\ell=\]\Plus_{r=0}^{2\ell}\C\)v_r$ as follows:
\vvn-.3>
$$
E\>v_r\>=\>\qnum r\)v_{r-1}\,,\qquad F\)v_r\>=\>\qnum{2\ell-r}\)v_{r+1}\,,
\qquad H\)v_r\>=\>(\ell-r)\>v_r\,.\kern-1.8em
\vv-.1>
\Tag{Vl}
$$
It is \irr/, if $q^{\)2k}\}\ne 1$ for any $k=1\lc 2\ell$. Recall that for such
\vvgood
$q$ the algebra $\Uq$ has a unique up to equivalence \irrep/ of dimension
$2\ell+1$. In what follows it is important for us that formulae \(Vl) are
analytic in $q$, though their explicit form is not quite essential.
\vsk.1>
Entries of the mondromy matrix
\ifMag\vvn-.6>\else\vvn-.2>\fi
$$
z_1\]\ldots\)z_N\>L_N(u\)/\]z_N)\)\ldots\>L_1(u\)/\]z_1)\,=\,
\pmatrix \)A(u) & B(u)\) \\\nn4> \)C(u) & D(u) \endpmatrix
\ifMag\vv-.3>\fi
$$
are \pol/s in $u$ and Laurent \pol/s in $q^{\)1/2}\}$ with values in
$\End(\Vox\])$. If necessary, their dependence on $q$ will be shown explicitly,
for instance, $A(u\);q)$.
\vsk.1>
The transfer-matrix \,$T(u)\>=\>A(u)\)+\)\ka\)D(u)$ commutes with the operator
$\Hs=H_1\]\lsym+H_N$ and is a generating \fn/ of the standard commuting
conservation laws. The algebraic \Ba/ gives a way to find \egv/s and \eva/s
of the transfer-matrix, \cf. Proposition~\[BA]. Recall the main points of
the technique.
\vsk.1>
The algebraic \Ba/ is based on the following commutation relations for
entries of the monodromy matrix:
\vvn-.7>
$$
\<B>
\gather
\bigl[\)B(u)\>,B(v)\)\bigr]\,=\,0\,,
\Tag{B}
\\
\nn7>
{\align
(u-v)\>A(u)\)B(v)\, &{}=
\,(u\)\q\}-v\)q)\>B(v)\)A(u)\>+\>v\>(q-\q)\>B(u)\)A(v)\,,
\ifMag\kern-1.3em\fi
\Tag{AB}
\\
\nn6>
(u-v)\>D(u)\)B(v)\, &{}=
\,(u\)q-v\)\q)\>B(v)\)D(u)\>-\>v\>(q-\q)\>B(u)\)D(v)\,,
\ifMag\kern-1.3em\fi
\Tag{DB}
\endalign}
\\
\ifMag\cnn-.1>\fi
\endgather
$$
and formulae for the action of $A(u)$ and $D(u)$ on the vector $\vox$:
\vvn-.1>
\ifMag
$$
\align
A(u)\> &\vox\,=\,\Ac(u)\>\vox\,,
\Tag{AD}
\\
\nn4>
& \Ac(u)\,=\,\prod_{i=1}^N\,(u\)q^{\)\ell_i}-z_i\)q^{-\ell_i})\,,
\\
\nn8>
D(u)\> & \vox\,=\,\Dc(u)\>\vox\,,\kern-1.6em
\Tag{DA}
\\
\nn4>
& \Dc(u)\,=\,\prod_{i=1}^N\,(u\)q^{-\ell_i}-z_i\)q^{\)\ell_i})\,.
\endalign
$$
\vsk-.4>
\else
$$
\gather
A(u)\>\vox\,=\,\Ac(u)\>\vox\,,\qqq D(u)\>\vox\,=\,\Dc(u)\>\vox\,,\kern-1.6em
\Tag{AD}
\\
\nn4>
\Ac(u)\,=\,\prod_{i=1}^N\,(u\)q^{\)\ell_i}-z_i\)q^{-\ell_i})\,,\kern 7em
\Dc(u)\,=\,\prod_{i=1}^N\,(u\)q^{-\ell_i}-z_i\)q^{\)\ell_i})\,.\kern-1.6em
\Tag{DA}\endgather
$$
\vsk-.2>
\fi
Let $\bv\tk\)=\>B(t_1)\ldots\)B(t_k)\>\vox$. Formulae \(B)\)--\)\(DA) imply
that the vector $\bv\tk$ is a \sym/ \fn/ of $\tk$, while the operators $\Hs$,
$A(u)$ and $D(u)$ act on it as follows:
\ifMag
\vvn-.1>
$$
\align
&\kern2.3em \Hs\>\bv\tk\,=\,(\ellsum\]-k)\,\bv\tk\,,
\\
\nn8>
A(u)\> & \bv\tk\,=\,
\Ac(u)\)\prod_{a=1}^k{u\)\q\}-t_a\)q\over u-t_a}\,\,\bv\tk\,+{}
\Tagg{A}
\\
\nn2>
{}+{}\, & (q-\q)\,\sum_{a=1}^k\>{t_a\over u-t_a\]}\;\Ac(t_a)
\prod_{\tsize{b=1\atop b\ne a}}^k{t_a\)\q\}-t_b\>q\over t_a\]-t_b}
\,\,\bv{u\),\tka}\,,\kern-.5em
\\
\nngood
\nn4>
\nngood
D(u)\> & \bv\tk\,=\,
\Dc(u)\)\prod_{a=1}^k{u\)q-t_a\)\q\}\over u-t_a}\,\,\bv\tk\,-{}
\Tagg{D}
\\
\nn2>
{}-\,{} & (q-\q)\,\sum_{a=1}^k\>{t_a\over u-t_a\]}\;\Dc(t_a)
\prod_{\tsize{b=1\atop b\ne a}}^k{t_a\)q-t_b\>\q\}\over t_a\]-t_b}
\,\,\bv{u\),\tka}\,.\kern-.5em
\endalign
$$
\vsk-.2>
\else
$$
\gather
\Hs\>\bv\tk\,=\,(\ellsum\]-k)\,\bv\tk\,,
\\
\nn8>
{\align
A(u)\>\bv\tk\, &{}=\,
\Ac(u)\)\prod_{a=1}^k{u\)\q\}-t_a\)q\over u-t_a}\,\,\bv\tk\,+{}
\Tag{A}
\\
\nn2>
&\){}+\,(q-\q)\,\sum_{a=1}^k\>{t_a\over u-t_a\]}\;\Ac(t_a)
\prod_{\tsize{b=1\atop b\ne a}}^k{t_a\)\q\}-t_b\>q\over t_a\]-t_b}
\,\,\bv{u\),\tka}\,,
\\
\nn2>
D(u)\>\bv\tk\, &{}=\,
\Dc(u)\)\prod_{a=1}^k{u\)q-t_a\)\q\}\over u-t_a}\,\,\bv\tk\,+{}
\Tagg{D}
\\
\nn2>
&\){}-\,(q-\q)\,\sum_{a=1}^k\>{t_a\over u-t_a\]}\;\Dc(t_a)
\prod_{\tsize{b=1\atop b\ne a}}^k{t_a\)q-t_b\>\q\}\over t_a\]-t_b}
\,\,\bv{u\),\tka}\,.
\endalign}
\endgather
$$
\vsk-.3>
\fi
The system of \Bae/s for rapidities $\tk$ has the form:
\vvn-.2>
$$
\Ac(t_a)\)\prod_{\tsize{b=1\atop b\ne a}}^k\,(t_a\]-t_b\>\qq)\,=\,
\ka\,\Dc(t_a)\)\prod_{\tsize{b=1\atop b\ne a}}^k\,(t_a\)\qq\}-t_b)\,,
\qqq\ a=1\lc k\,.\kern-4em
\ifMag\vv-.3>\else\vv-.4>\fi
\Tag{xxz}
$$
We do not distinguish \sol/s of this system which are obtained from each other
by \perm/s of the \var/s $\tk$. We say that a \sol/ $\tk$ contains a point $u$
if $u\in\lb\)\tk\rb$. A \sol/ $\tk$ of system \(xxz) is called \em{offdiagonal}
if $t_a\]\ne t_b$ for all $a\),b=1\lc k$. The vector $\bv\tk$ is called
the \Bv/ if $\tk$ is an offdiagonal \sol/ of the \Bae/s. The vector $\vox$
is the \Bv/ corresponding to the empty set of rapidities by convention.
\Prop{BA}
Let $\tk$ be an offdiagonal \sol/ of system \(xxz)\). Then
\vvn-.1>
\ifMag
$$
\align
& T(u)\>\bv\tk\,={}
\Tag{Tbv}
\\
\nn4>
& \,{}=\,\Bigl(\)\Ac(u)\)\prod_{a=1}^k{u\)\q\}-t_a\)q\over u-t_a}
\;+\,\ka\>\Dc(u)\)\prod_{a=1}^k{u\)q-t_a\)\q\}\over u-t_a}\,\Bigr)\,
\bv\tk\,.\kern-1.9em
\endalign
$$
\else
$$
T(u)\>\bv\tk\,=\,\Bigl(\)\Ac(u)\)\prod_{a=1}^k{u\)\q\}-t_a\)q\over u-t_a}
\;+\,\ka\>\Dc(u)\)\prod_{a=1}^k{u\)q-t_a\)\q\}\over u-t_a}\,\Bigr)\,
\bv\tk\,.
\Tag{Tbv}
$$
\fi
\endpro
\nt
The statement follows from formulae \(A) and \(D)\).
\Par
Henceforth we assume that $q^{\)2\ell_i}\]z_i\ne q^{-2\ell_j}\]z_j$ for
any $i\),j=1\lc N$, that is, the \pol/s $\Ac(u)$ and $\Dc(u)$ are coprime.
Furthermore, we assume that $q^{\)2r}\}\ne 1$ for all
$r=1\lc\max\)(2\)\ell_1\lc 2\)\ell_N\])$, in other words,
that all the \rep/s $V_{\ell_1}\lc V_{\ell_N}$ are \irr/.
\vsk.1>
A \sol/ $\tk$ of system \(xxz) is called \em{admissible} if $t_a\]\ne 0$ and
$t_a\]\ne\qq t_b$ for all $a\),b=1\lc k$. Simple analysis of system \(xxz)
yields the following properties of admissible \sol/s
\Lm{2}
Any admissible \sol/ of system \(xxz) does not contain the points $\zlpmn$.
\endpro
\Lm{3}
Let $q^{\)2m}\}\ne 1$ for all \>$m=1\lc k$.
\atem
If a \sol/ of system \(xxz) does not contain any of the points $\zlpn$,
then it is admissible.
\bitem
If a \sol/ of system \(xxz) does not contain any of the points $\zlmn$,
then it is admissible.
\endpro
We say that the points \)$\zn$ are \em{well separated},
if $q^{\)2(r-\)\ell_i)}\]z_i\ne q^{\)2(s-l_j)}\]z_j$
for all $r=0\lc 2\ell_i$, $s=0\lc 2\ell_j$, $i\),j=1\lc N$.
\Th{TV}
\back\cite{{\rm TV}, Theorems~4.2 \)and 5.1\)}
Let the points \)$\zn$ be well separ\-ated. Then for generic $\ka$ all
admissible offdiagonal \sol/s of system \(xxz) are \ndeg/, the number of
distinct admissible offdiagonal \sol/s of system \(xxz) equals the dimension
of the subspace $\bigl\lb\)v\in\Vox\vert{}\Hs\)v\)=(\ellsum\]-k)\>v\)\bigr\rb$,
and the corresponding \Bv/s form a basis of this subspace. If $q$ is not a root
of unity, then the \Bv/s corresponding to inadmissible offdiagonal \sol/s of
system \(xxz) equal zero.
\endpro
It is plausible that the assumption of Theorem~\[TV] about the points $\zn$
being well separated can be weakened. For example, the statement of
the theorem conjecturally remains true for $z_1\]\lsym=z_N$.

\Shead
Assume that $\qq\}$ is an \}\$M$-th root of unity. In principle, for $k\ge M$
system \(xxz) can have inadmissible \sol/s of the form
$\cyc{u\)},\)t_{M+1}\lc t_k$ with arbitrary $u$; in particular, such \sol/s
are not isolated. The \Bv/ corresponding to such a \sol/ equals zero, because
$$
B(u)\>B(u\)\qq)\)\ldots\>B(u\)q^{\)2M-\)2})\,=\,0
$$
for any $u$, see \cite{T}\). Besides, the sequence $t_{M+1}\lc t_k$ is a \sol/
of system \(xxz) for $\)k-\]M$ \var/s.
\vsk.1>
Theorem~\[TV] indicates that to find the spectrum of the transfer-matrix $T(u)$
for generic $\ka$ it suffices in general to consider only \Bv/s corresponding
to admissible \sol/s of system \(xxz)\). At the same time, one can see from
the results of \cite{FM2} that for special values of $\ka$ inadmissible \sol/s
of system \(xxz) mentioned above and corresponding them analogues of \Bv/s can
play essential role in constructing \egv/s of the transfer-matrix.
\vsk.1>
Assume that $\ka=q^{\)2(p\)-\)\ellsum)}\}$ for certain integer $p$. Let $\tk$
be an offdiagonal \sol/ of system \(xxz)\). The main aim of this note is to
construct vectors $\bbv{\tk\);\)\um}_p$, depending on the parameters $\um$
and such that
\ifMag
$$
\alignat2
&&\Llap{\Hs\,\bbv{\tk\);\)\um}_p\,=\,(\ellsum\]-k-mM)\,\bbv{\tk\);\)\um}_p\,,}&
\\
\nn14>
T &{}(u)\>\bbv{\tk\);\)\um}_p\,={}
\Tagg{Tbbv}
\\
\nn4>
{}=\,q^{\)mM}\, &\}\Bigl(\)\Ac(u)\)\prod_{a=1}^k{u\)\q\}-t_a\)q\over u-t_a}
\;+\,\ka\>\Dc(u)\)\prod_{a=1}^k{u\)q-t_a\)\q\}\over u-t_a}\,\Bigr)\,
\bbv{\tk\);\)\um}_p\,,
\endalignat
$$
\else
$$
\align
\Hs\, & {}\bbv{\tk\);\)\um}_p\,=\,(\ellsum\]-k-mM)\,\bbv{\tk\);\)\um}_p\,,
\\
\nn10>
T(u)\> & {}\bbv{\tk\);\)\um}_p\,={}
\Tag{Tbbv}
\\
\nn4>
&\]\!\!{}=\,q^{\)mM}\)\Bigl(\)\Ac(u)\)\prod_{a=1}^k{u\)\q\}-t_a\)q\over u-t_a}
\;+\,\ka\>\Dc(u)\)\prod_{a=1}^k{u\)q-t_a\)\q\}\over u-t_a}\,\Bigr)\,
\bbv{\tk\);\)\um}_p\,,\kern-1.8em
\\
\cnn-.2>
\endalign
$$
\fi
\cf. Proposition~\[BAM]. Though the parameters $\um$ are arbitrary, this does
not contradict to the \fd/ity of the space of states ${\Vox}$, since the \eva/
of $T(u)$ does not depend on $\um$. The construction being suggested
generalizes that of \cite{FM2}\).
\vsk.1>
The vectors $\bbv{\tk;\um}_p$ with different $m$ correspond to the same \eva/
of the operators \>$q^{\)\pm\Hs}\)T(u)$. The respective \eva/s of the
transfer-matrix $T(u)$ coincide, if $q^{\)mM}\}=1$, and can differ by the sign,
if $q^{\)mM}\}=-1$. Comparison formulae \(Tbv) and \(Tbbv) shows that
the vector $\bbv{\tk;\um}_p$ can be viewed as an analogue of the \Bv/
corresponding to the inadmissible \sol/
$\tk,\)u_1\)\lc\)u_1\)q^{\)2M-\)2}\lc\)u_m\)\lc\)u_m\)q^{\)2M-\)2}\}$
of system \(xxz) for $k+mM$ \var/s.

\Shead
Fix an integer ${M>1}$, and let ${\gm_0=\pi K\]/M}$ for certain $K$ coprime
with $M$. Set $\qo=e^{i\gm_0}\}$ and $\eta=q\)/\]\qo$. For any object depending
on $q$ we assume that $q=\qo$, unless the dependence is shown explicitly.
\vsk.2>
Taking into account that
$B(u)\>B(u\)\qqz)\)\ldots\>B(u\)\qz^{\)2M-\)2})\,=\,0$, introduce
\vvn.1>
an operator $\BB\)(u\);X)$ depending on a vector $X=(\xm)\in\C^{\)M}\}$:
\vvn.2>
$$
\BB\)(u\);X)\,=\lim_{q\)\to\)\qo\!}\bigl((q-\qo)\1 B(u\)\eta^{x_1};q)
\>B(u\)\qq\eta^{x_2};q)\)\ldots\>B(u\)q^{\)2M-\)2}\eta^{x_M};q)\bigr)\,,
\kern-.8em
\Tag{limB}
$$
assuming that $\eta^x\}\to 1$ as $\eta\to 1$. The operator $\BB\)(u\);X)$
does not change under the simultaneous shift of all the parameters $\xm$ by
the same number. Explicit calculation of the limit in formula \(limB) yields
\ifMag
$$
\align
\BB\)(u\);X)\,=\)\sum_{r=0}^{M-1}B(u)\ldots\)B(u\)\qz^{\)2r-2})\)
\bigl(\)\der_q B(u\)\qz^{\)2r})+
u\)\qz^{\)2r-1}(x_{r+1}\]+2r)\)\der_u B(u\)\qz^{\)2r})\bigr)\)\x{}\!\} &
\\
\nn3>
{}\x\)B(u\)\qz^{\)2r+2})\ldots\)B(u\)\qz^{\)2M-\)2})\,, &
\endalign
$$
\else
\vvn.3>
$$
\Lline{\BB\)(u\);X)\,=\)\sum_{r=0}^{M-1}B(u)\ldots\)B(u\)\qz^{\)2r-2})\)
\bigl(\)\der_q B(u\)\qz^{\)2r})+
u\)\qz^{\)2r-1}(x_{r+1}\]+2r)\)\der_u B(u\)\qz^{\)2r})\bigr)
\)B(u\)\qz^{\)2r+2})\ldots\)B(u\)\qz^{\)2M-\)2})\,,}
\vv.1>
$$
\fi
which is similar to formula (1.38) in \cite{FM2}\). However, it is much more
convenient to use formula \(limB) itself, taking the limit $q\to\qo\}$ only at
the end of computation. In particular, relation \(B) immediately implies that
\vvn.4>
$$
\bigl[\)\BB\)(u\);X)\>,B(v)\)\bigr]\,=\,
\bigl[\)\BB\)(u\);X)\>,\BB\)(v\);Y)\)\bigr]\,=\,0
\vv.3>
$$
for any $u\),v\),X\),Y$ and, therefore, the vector
\ifMag\else\vvn.3>\fi
$$
\bbbv{\tk;\um\);\Xm}\>=\,\BB\)(u_1\);X_1)\ldots\)\BB\)(u_m\);X_m)\>\bv\tk
\ifMag\else\vv.1>\fi
$$
is a \sym/ \fn/ of $\tk$, as well as the pairs $(u_1\),X_1)\lc (u_m\),X_m)$.
It is clear that
\ifMag\vvn-.1>\else\vvn-.7>\fi
$$
\align
\Hs\, &{} \bbbv{\tk;\um\);\Xm}\>={}
\\
\nn8>
& {}\!\!=\,(\ellsum\]-k-mM)\,\bbbv{\tk;\um\);\Xm}\>.
\endalign
$$
\par
Consider the \pol/ \;$\dsize P(u)\,=\,
\prod_{i=1}^N\,\prod_{r=0}^{2\ell_i-1}(u-z_i\)\qz^{\)2(\ell_i-\)r)})$
\;and the \fn/s
\vvn-.5>
$$
\gather
Q_n(u\);\tk)\,=\;
{u^{\)k-n}\)P(u)\over\prod_{a=1}^k\)(u-t_a)\>(u-t_a\)\qqz)}\;,
\\
\nn8>
F_n(u\);\tk)\,=\;{1\over M}\,\sum_{r=0}^{M-1} Q_n(u\)\qz^{\)2r};\tk)\,,
\\
\nn7>
G_n(u\);\tk)\,=\;{1\over M}\,\sum_{r=1}^{M-1}\>r\>Q_n(u\)\qz^{\)2r};\tk)\,,
\\
\cnn.4>
\nngood
\cnn-.5>
\endgather
$$
which satisfy the relations
\vvn-.5>
$$
\gather
{Q_n(u\)\qqz)\over Q_n(u)}\;=\,\qz^{\)2(\ellsum\)-\>n)}\,\,
{\Ac(u)\over \Dc(u)}\;\prod_{a=1}^k\>{u-t_a\)\qqz\over u\)\qqz-t_a}\;,
\\
\ifMag\nn10>\else\nn9>\fi
F_n(u\)\qqz)\,=\,F_n(u)\,,\qqq
G_n(u\)\qqz)\)-\)G_n(u)\,=\,Q_n(u)\)-\)F_n(u)\,.
\\
\cnn.1>
\endgather
$$
If $\tk$ is an admissible offdiagonal \sol/ of system \(xxz), then it is easy
to see that the \fn/ $F_n(u\);\tk)$ is a \pol/ in $u$.
\vsk.1>
Set \,${{\bbv{\tk;\um}_n\,=\,\bbbv{\tk;\um\);X_1\"n\lc X_m\"n}}}$ where
\vvn.06>
the vectors $X_i\"n\]=\>X\"n(u_i\);\tk)\>=\>(x_1\"n\}\lc x_M\"n)(u_i\);\tk)$,
\;$i=1\lc m$, are defined as follows:
$$
x_r\"n(u\);\tk)\,=\,2\)\bigl(1-r-G_n(u\)q^{\)2r-2};\tk)/F_n(u\);\tk)\bigr)\,.
\kern-2em
\vv.2>
\Tag{xir}
$$
\Prop{BAM}
Let ${\ka=\qz^{\)2(p\)-\)\ellsum)}}\}$ for a certain integer $p$, and let
$\tk$ be an offdiagonal \sol/ of system \(xxz) at ${q=\qo}$. Assume that
\,$F_p(u_i\);\tk)\ne 0$ \>and \,$u_i^M\!\ne t_a^M$ \)for all \>$a=1\lc k$,
$i=1\lc m$. Then
\vvn.1>
\ifMag
$$
\align
& T(u)\>\bbv{\tk\);\)\um}_p\,={}
\\
\nn5>
&\}\!{}=\,\qz^{\)mM}\)\Bigl(\)\Ac(u)\)
\prod_{a=1}^k{u\)\qz\1\}-t_a\)\qo\over u-t_a}
\;+\,\ka\>\Dc(u)\)\prod_{a=1}^k{u\)\qo-t_a\)\qz\1\}\over u-t_a}\,\Bigr)\,
\bbv{\tk\);\)\um}_p\,.
\\
\cnn-.4>
\endalign
$$
\else
$$
\align
T(u)\> & {}\bbv{\tk\);\)\um}_p\,={}
\\
\nn5>
&\]\!\!{}=\,\qz^{\)mM}\)\Bigl(\)\Ac(u)\)
\prod_{a=1}^k{u\)\qz\1\}-t_a\)\qo\over u-t_a}
\;+\,\ka\>\Dc(u)\)\prod_{a=1}^k{u\)\qo-t_a\)\qz\1\}\over u-t_a}\,\Bigr)\,
\bbv{\tk\);\)\um}_p\,.\kern-2em
\\
\cnn-.4>
\endalign
$$
\fi
\endpro
\Pf.
Using formulae \(A) and \(D) for generic $q$ and taking the limit $q\to\qo$
we obtain the following formulae for the action of the operators $A(u)$ and
$D(u)$ on the vector $\bbbv{\tk;\um\);\Xm}$:
$$
\alignat2
\ifMag\else A(u)\>\fi &{} \ifMag \]A(u)\>\fi \bbbv{\tk;\um\);\Xm}\,={}
\Tag{A2}
\\
\nn7>
& {}\!\!=\,\qz^{\)mM}\>\biggl[\,\Ac(u)\)
\prod_{a=1}^k{u\)\qz\1\}-t_a\)\qo\over u-t_a}\,\,\bbbv{\tk;\um\);\Xm}\,+{}
\\
\nn3>
& {}\!\}+\,(\qo-\qz\1)\>\sum_{a=1}^k\,{t_a\over u-t_a\]}\;\Ac(t_a)
\prod_{\tsize{b=1\atop b\ne a}}^k{t_a\)\qz\1\}-t_b\>\qo\over t_a\]-t_b}\;\x{}
\\
\nn7>
&& \Llap{{}\x\,\bbbv{u\),\tka;\um\);\Xm}\,-{}} &
\\
\nn7>
& {}\!\}-\>\sum_{i=1}^m\,\sum_{r=0}^{M-1}
{u_i\)\qz^{\)2r}\over u-u_i\)\qz^{\)2r}\]}\;
(x_{i,r+1}\]-x_{ir})\,\Ac(u_i\)\qz^{\)2r})
\prod_{a=1}^k{u_i\)\qz^{\)2r-1}\}-t_a\)\qo\over u_i\)\qz^{\)2r}\}-t_a}\;\x{}
\\
\nn4>
& \Rlap{\Lph{{}\!\]-\>\sum_{i=1}^m}{{}\x{}}\,
\prod_{\tsize{s=0\atop s\ne r}}^{M-1}B(u_i\)\qz^{\)2s})\,
\bbbv{u\),\tk;\umi\);\Xmi}\>\biggr]\,,}
\\
\cnn-.7>
\nngood
\cnn.7>
\endalignat
$$
\vsk-.6>\nt
$$
\alignat2
\ifMag\else D(u)\>\fi &{} \ifMag \}D(u)\>\fi \bbbv{\tk;\um\);\Xm}\,={}
\Tag{D2}
\\
\nn7>
& {}\!\!=\,\qz^{\)mM}\>\biggl[\,\Dc(u)\)
\prod_{a=1}^k{u\)\qo-t_a\)\qz\1\over u-t_a}\,\,\bbbv{\tk;\um\);\Xm}\,-{}
\\
\nn3>
& {}\!\}-\,(\qo-\qz\1)\>\sum_{a=1}^k\,{t_a\over u-t_a\]}\;\Dc(t_a)
\prod_{\tsize{b=1\atop b\ne a}}^k{t_a\)\qo-t_b\>\qz\1\over t_a\]-t_b}\;\x{}
\\
\nn7>
&& \Llap{{}\x\,\bbbv{u\),\tka;\um\);\Xm}\,+{}} &
\\
\nn7>
& {}\!\}+\>\sum_{i=1}^m\,\sum_{r=0}^{M-1}
{u_i\)\qz^{\)2r}\over u-u_i\)\qz^{\)2r}\]}\;
(x_{i,r+2}\]-x_{i,r+1})\,\Dc(u_i\)\qz^{\)2r})
\prod_{a=1}^k{u_i\)\qz^{\)2r+1}\}-t_a\)\qz\1\over u_i\)\qz^{\)2r}\}-t_a}\;\x{}
\\
\nn4>
& \Rlap{\Lph{{}\!\]+\>\sum_{i=1}^m}{{}\x{}}\,
\prod_{\tsize{s=0\atop s\ne r}}^{M-1}B(u_i\)\qz^{\)2s})\,
\bbbv{u\),\tk;\umi\);\Xmi}\>\biggr]\,.}
\\
\cnn-.4>
\endalignat
$$
Here \,$X_i=(\xim)$, \,$x_{i0}=x_{iM}\]+2M$, \,$x_{i,M+1}=x_{i1}\]-2M$.
\vsk.3>
The action of the transfer-matrix $T(u)$ on the vector $\bbbv{\tk;\um\);\Xm}$
produces ``unwanted terms'' of two kinds, which are \resp/ generated by
the second (ordinary sums) and third (double sums) terms in \rhs/s of formulae
\(A2) and \(D2). The unwanted terms of the first kind cancel each other
if $\tk$ is an offdiagonal \sol/ of system \(xxz) at $q=\qo$. Regardless of
the cancelation of the unwanted terms of the first kind, the unwanted terms of
the second kind cancel each other if
\ifMag
\vv-.4>
$$
\align
(x_{i,r+1}\] &{} -x_{ir})\,\Ac(u_i\)\qz^{\)2r})
\prod_{a=1}^k{u_i\)\qz^{\)2r}\}-t_a\)\qqz\over u_i\)\qz^{\)2r}\}-t_a}\;={}
\Tag{xeq}
\\
& {}=\,\ka\,(x_{i,r+2}\]-x_{i,r+1})\,\Dc(u_i\)\qz^{\)2r})
\prod_{a=1}^k{u_i\)\qz^{\)2r+2}\}-t_a\over u_i\)\qz^{\)2r}\}-t_a}
\endalign
$$
\else
$$
\kern-.4em
(x_{i,r+1}\]-x_{ir})\,\Ac(u_i\)\qz^{\)2r})
\prod_{a=1}^k{u_i\)\qz^{\)2r}\}-t_a\)\qqz\over u_i\)\qz^{\)2r}\}-t_a}\;=\,
\ka\,(x_{i,r+2}\]-x_{i,r+1})\,\Dc(u_i\)\qz^{\)2r})
\prod_{a=1}^k{u_i\)\qz^{\)2r+2}\}-t_a\over u_i\)\qz^{\)2r}\}-t_a}
\Tag{xeq}
$$
\fi
for all $r=0\lc M-1$, $i=1\lc m$. For $\ka=q^{\)2(p\)-\)\ellsum)}$
these \eq/s hold if $x_{ir}=\)x_r\"p(u_i)$ for all $r=1\lc M$, $i=1\lc m$.
\epf
For the homogenuous spin\>-$\!{1\over 2}$ \XXZ/ model (six-vertex model) with
the even number of lattice vertices and the periodic boundary conditions:
$\ell_1\]\lsym=\ell_N=1/2$, $z_1\]\lsym=z_N=1$, $\ka=1$, Proposition~\[BAM]
was proved in \cite{FM2}\).
\Rem
Consider relations \(xeq) for a fixed $i$ as \eq/s for $\xim$ with
the boundary conditions $x_{i0}=x_{iM}\]+2M$, $x_{i,M+1}=x_{i1}\]-2M$,
\vvn.06>
the \var/s $u_i\),\,\tk$ being given. Assume that $u_i^M\!\ne 1$ and
$t_a^M\!\ne u_i^M\}$ for all \>$a=1\lc k$. It is easy to check that
under these assumptions system \(xeq) can have \sol/s only if
$\ka^M\!=\qz^{\)2M(\ellsum)}\}$, and one has $x_{ir}\ne x_{i,r+1}$
for all $r=0\lc M$.
\vsk.1>
Denote
\vvn.2>
$$
y_r\,=\;{\Ac(u_i\)\qz^{\)2r-2})\over\ka\>\Dc(u_i\)\qz^{\)2r-2})}\;
\prod_{a=1}^k\>{u_i\)\)\qz^{\)2r-2}\}-t_a\)\qqz\over u_i\)\qz^{\)2r}\}-t_a}\;.
\ifMag\vv.2>\fi
$$
Then \,$x_{i,r+1}\]-x_{ir}\)=\)y_1\>y_2\]\ldots\)y_r\>(x_{i1}-x_{i0})$ \,and
\ifMag
\vvn-.1>
$$
(x_{i1}-x_{i0})\tsum_{r=0}^{M-1}\]y_1\>y_2\]\ldots\)y_r\)=\)
x_{iM}-x_{i0}\)=\)-\)2\)M\,,
\vv-.2>
$$
\else
\>$(x_{i1}-x_{i0})\sum_{r=0}^{M-1}\]y_1\>y_2\]\ldots\)y_r\)=\)
x_{iM}-x_{i0}\)=\)-\)2\)M$,
\fi
which implies that the general \sol/ of the system in question has the form
$$
x_{ir}\>=\,x_{i1}\)-\>2\)M\)\tsum_{s=1}^{r-1} y_1\>y_2\]\ldots\)y_s\,
\Bigl(\>\tsum_{s=0}^{M-1}\]y_1\>y_2\]\ldots\)y_s\Bigr)^{\!-1},\qqq r=2\lc M\,,
\kern-2em
\vv-.3>
$$
with arbitrary $x_{i1}$. If $\sum_{r=0}^{M-1}\]y_1\>y_2\]\ldots\)y_r\)=\)0$,
then system \(xeq) is not solvable.
\vsk.1>
Thus, under the natural assumptions the solvability of system \(xeq)
implies that $\ka=\qz^{\)2(p\)-\)\ellsum)}\}$ for a certain integer $p$.
In this case the \sol/ of system \(xeq) has the form
\vvn-.4>
$$
x_{ir}\>=\>c_i+\)x_r\"p(u_i\);\tk)\,, \qqq r=1\lc M\,,\qquad i=1\lc m\,,
\kern-4em
\vv-.1>
$$
\cf. \(xir)\), with arbitrary $\cm$, while the vector $\bbbv{\tk;\um\);\Xm}$
does not depend on $\cm$ and equals $\bbv{\tk;\um}_p$.
\enddemo
\Rem
It has been shown in \cite{DFM} that there is an action of the loop algebra
$\slt$ in the space of states of the six-vertex model at roots of unity,
and this action commutes with the transfer-matrix of the six-vertex model.
The existence of a large symmetry algebra causes the degeneration of the
spectrum of the transfer-matrix. As mentioned in \cite{FM2} the symmetry
degeneration of the spectrum of the transfer-matrix of the six-vertex model
apparently corresponds to that degeneration of the spectrum which occurs in
this case due to Proposition~\[BAM]. Since Proposition~\[BAM] remains true and
for the inhomogeneous arbitrary spin \XXZ/ model with suitable quasiperiodic
boundary conditions, one may suppose that the model has a large symmetry
algebra in the general case too. Moreover, it could happen that an action of
this algebra exists for any quasiperiodic boundary conditions, but commutes
with the transfer-matrix only for special values of the quasiperiodicity
parameter.

\enddemo

\vsk>
\myRefs
\widest{VVV}
\parskip.1\bls

\ref\Key DFM
\by T\]\&Deguchi, K\&Fabricius and B\&M\&McCoy
\paper The $sl_2$ loop algebra symmetry of the six vertex model at roots of
unity \jour J\.Stat\.Phys. \vol 102 \yr 2001 \issue 3\)\~\)4 \pages 701\~\)736
\endref

\ref\Key FM1
\by K\&Fabricius and B\&M\&McCoy
\paper Bethe's equation is incomplete for the XXZ model at roots of unity
\jour J\.Stat\.Phys. \vol 103 \yr 2001 \issue 5\)\~\>6 \pages 647\)\~\>678
\endref

\ref
\by K\&Fabricius and B\&M\&McCoy
\paper Completing Bethe's equation at roots of unity
\jour J\.Stat\.Phys. \vol 104 \yr 2001 \issue 3\)\~\)4 \pages 573\)\~\)587
\endref

\ref\Key FM2
\by K\&Fabricius and B\&M\&McCoy
\paper Evaluation parameters and Bethe roots for the six vertex model at roots
of unity \jour Preprint YITPSB\)\~\>01\~\>42 \yr 2001 \pages 1\~\)24
\endref

\ref\Key KBI
\by \Kor/, N\&M\&Bogoliubov and \Izergin/
\book Quantum inverse scattering method and correlation \fn/s
\yr 1993 \publ \CUP/ \pages 556\,
\endref

\goodsk:->

\ref\Key T
\by \VT/
\paper Cyclic monodromy matrices for the \Rm/ of the \XXZ/-model
and the chiral Potts model with fixed-spin boundary conditions
\jour Int\.J\.Mod\.Phys\. A (\rm Suppl\.1B) \yr 1992 \vol 7 \pages 963\)\~\)975
\endref

\ref
\by \VT/
\paper Cyclic monodromy matrices for $sl(n)$ trigonometric \Rms/
\jour \hbox{\CMP/} \vol 158 \yr 1993 \issue 3 \pages 459\>\~\)483
\endref

\ref\Key TV
\by \VT/ and \Varch/
\paper Completeness of \Bv/s and \deq/s with regular singular points
\jour \IMRN/ \yr 1996 \issue 13 \pages 637\)\~\>669
\endref

\endRefs

\bye